\documentclass[final,5p,times,twocolumn,authoryear]{elsarticle}

\usepackage{amssymb}

\usepackage{amsmath}
\usepackage{tikz}
\usepackage{subcaption}
\usepackage{xcolor}
\usepackage{mathtools}
\usepackage{layouts}
\usepackage{enumitem}

\usepackage{amsmath,amssymb,amsfonts}
\usepackage{algorithmic}
\usepackage{graphicx}
\usepackage{textcomp}
\usepackage{array}
\usepackage{booktabs}
\usepackage{multirow}

\usepackage{url}

\usepackage{breakurl}
\usepackage[breaklinks]{hyperref}

\newenvironment{newtext}
  {\begingroup\color{black}}
  {\endgroup}

\usepackage{nomencl}
\makenomenclature

\usepackage{etoolbox}
\renewcommand\nomgroup[1]{%
  \item[\bfseries
  \ifstrequal{#1}{P}{Parameters}{%
  \ifstrequal{#1}{V}{Variables$^*$}{%
    \ifstrequal{#1}{S}{Sets}{%
        \ifstrequal{#1}{F}{Functions}{%
  \ifstrequal{#1}{O}{Other symbols}{}}}}}%
]}

\newcommand{\nomunit}[1]{%
\renewcommand{\nomentryend}{\hspace*{\fill}#1}}

\journal{Elsevier}

\begin{document}

\begin{frontmatter}



\title{The Value of Ancillary Services for Electrolyzers}


\author[label1]{Andrea Gloppen Johnsen}
\author[label1]{Lesia Mitridati}
\author[label2]{Donato Zarrilli}
\author[label1]{Jalal Kazempour}

\affiliation[label1]{organization={Department of Wind and Energy Systems, 
 Technical University of Denmark},
            addressline={}, 
            city={ Kgs. Lyngby},
            country={Denmark}}

\affiliation[label2]{organization={DNV Energy Systems},
            city={Bologna},
            country={Italy}}

\begin{abstract}
Although primarily designed for hydrogen production, electrolyzers can support power systems by providing various ancillary services, opening new revenue streams that enhance their economic viability.  This paper investigates the participation of an electrolyzer in frequency-supporting reserve markets, analyzing how bid structures and activation intensities affect its value. We develop a mixed-integer linear program to co-optimize electricity procurement and reserve provision, and analytically derive the opportunity cost of reserve provision, which determines the optimal bid price. Using historical price and frequency data from western Denmark, we show that asymmetric, hourly reserve products often entail no opportunity cost and can increase profits by up to 47\%. However, energy-intensive reserves may disrupt hydrogen production and risk unmet demand. Our findings reveal that flexible bidding can mitigate these risks while maintaining profitability. We also highlight the benefits of diversifying across reserve products and offer two recommendations: System operators should reconsider reserve bid structures to better accommodate electrolyzers, and electrolyzer owners should not overlook energy-intensive reserve services when hydrogen demand is flexible.

\end{abstract}




\end{frontmatter}


\vspace{2mm}
\section*{List of Abbreviations}

\begin{description}[labelwidth=1cm, labelsep=1cm]
    \item[DK1] Western Denmark Bidding Zone 
    \item[FCR] Frequency Containment Reserve 
    \item[mFRR] Manual Frequency Restoration Reserve  
\end{description}

\section{Introduction}


Electrolyzers and power-to-X facilities are increasingly acknowledged as key enablers of the green energy transition, providing viable pathways to indirectly electrify hard-to-abate sectors such as heavy transport and the chemical industry through the production of electrolytic hydrogen and e-fuels \citep{PtX_strategy_EU}. Beyond their decarbonization potential, electrolyzers could provide valuable flexibility in power systems dominated by stochastic renewable energy sources, where flexibility scarcity is a growing concern \citep{Anderson2004, Xiong2021}. 
However, uncertainty surrounding their long-term profitability remains a major barrier to their rapid deployment \citep{iea_hydrogen_review}. The provision of ancillary services,\footnote{The terms ``ancillary services" and ``reserves" are used interchangeably throughout this paper.} such as grid balancing reserves, has been identified as a promising solution to enhance their economic viability \citep{Kountouris2023, energinet2}.

The growing demand for ancillary services, combined with the declining share of thermal power generation as the main conventional source of power system flexibility, underscores the urgency of integrating ``new" technologies such as electrolyzers to supply up-regulation reserves \citep{energinet_AS_outlook}. A critical factor in their reserve market participation is the ability of electrolyzers to respond to frequency deviations by adjusting power consumption. While they can modify their power consumption set-points within ramping capabilities, such adjustments directly affect hydrogen production. For certain reserve products, frequent or prolonged activations can result in significant deviations in hydrogen output, posing operational challenges when production is constrained by downstream requirements such as purchase agreements, on-site storage limitations, or integration with processes like methanol production. Without accounting for activation impacts, reserve participation may lead to suboptimal or infeasible real-time operations. However, forecasting reserve activation, which is linked to the system imbalance and imbalance price, can be challenging at the time of bidding in the day-ahead stage \citep{Browell2022}.

Although ancillary services have the potential to improve profitability, electrolyzers are primarily designed to meet hydrogen demand rather than serve the grid for balancing purposes. This distinguishes them from other modern flexibility providers, such as batteries, which inherently operate as flexible resources with minimal net energy consumption over time. In contrast, electrolyzers must meet specific hydrogen production targets, which limits their ability to reduce energy consumption for grid services. However, electrolyzers may tolerate short-term variations in energy consumption more effectively than batteries.  Unlike batteries, electrolyzers may face contractual rather than physical constraints on energy consumption, allowing for some operational flexibility, albeit with potential monetary penalties. Recent requirements by Energinet, the Danish transmission system operator, for technologies with limited energy resources further highlight the need for energy-side flexibility in reserve provision \citep{energinet_tender_2024}.

This paper evaluates the value of ancillary services for an electrolyzer that serves a hydrogen demand, considering the impact of reserve activation on the final hydrogen schedule. The analysis focuses on frequency containment reserve (FCR) and manual frequency restoration reserve (mFRR) markets, which represent distinct reserve characteristics. FCR requires fast responses over short durations and is not energy-intensive, which makes it less impactful on hydrogen production. Conversely, mFRR involves slower but prolonged activations, significantly affecting hydrogen output. We conduct an ex-post analysis, utilizing historical market and frequency data from the western Denmark bid zone in the years 2021-2023, which provides an upper bound on electrolyzer profitability under observed market conditions. In addition, we analytically compare how the reserve bid structures impact the electrolyzer's willingness to participate in said reserve markets. These results provide insights into the value of ancillary services for electrolyzers, regardless of specific market conditions. 

The literature on the economic value of ancillary services for electrolyzers remains limited. For example, a preliminary analysis by Energinet \citep{Energinet} finds that participating in Danish reserve markets in 2021 could lower the levelized cost of hydrogen with $\sim$18\%. In a follow up report \citep{energinet2}, Energinet explores the additional value of having a flexible hydrogen off-take structure. They find that in 2023, offering reserves could lower the levelized cost of hydrogen from 6.20 to 5.15  EUR/kg, and that a flexible off-take structure yields a further reduction to 4.01 EUR/kg. However, both reports employ a simplified power-to-hydrogen conversion factor, neglecting electrolyzer-specific characteristics and the impact of reserve activation on hydrogen production.   In \cite{marco}, the value of FCR services for an electrolyzer in the Nordic synchronous region is estimated using a detailed physical model, demonstrating that FCR participation could account for up to 72\% of total revenues in past years with high FCR prices. However, this study does not consider reserve products beyond FCR, particularly those that are generally energy-intensive.  \cite{Gu2024} finds that the capacity payments from reserve provisions can exceed the revenues from hydrogen sales for a hydrogen production system participating in Nordic markets. Similarly, \cite{Pavi2022}, \cite{France}, \cite{Ennassiri2024}, and \cite{Germany_offshore} analyze the economic benefits of electrolyzers providing ancillary services in Croatia, France, Italy and Germany, respectively. Reserve activation is discussed in the literature, however, not in the context of energy-intensive reserve products. \cite{Pavi2022} address the activation of a reserve on a hydrogen network balance, ensuring feasibility of the network for maximum and minimum possible activations. \cite{Khajeh2024} investigates a wind-electrolyzer system that delivers FCR in Finland, ensuring that the reserved capacity can be fully activated under large deviations in real-time wind power output.  In contrast, \cite{ZHENG202326046} addresses reserve activation through chance constraints, ensuring feasibility of meeting the reserved capacity with high probability while allowing overbidding. This approach relies on accurate probabilistic forecasts of the activation stage and does not consider the impact of activation on meeting hydrogen demands. \cite{Cheng2023} investigates a joint wind and electrolyzer system providing frequency services and considers the effect of reserve activation on the final hydrogen production schedule. However, this work assumes that the activation of the reserve is inherently symmetric (i.e., no impact on the production schedule across time) and that the asymmetry in activation arises from the coordinated response of the wind and electrolyzer plants, and can thus be estimated accurately.   \cite{Naughton2023} models a virtual hybrid renewable hydrogen power plant that participates in the Australian energy and ancillary service markets. They embed a time duration factor of reserves to ensure that hydrogen storages are within feasible limits every half-hour time step. However, the reserve products considered have activation duration times of up to 5 minutes only.

The literature lacks investigations of the effects of reserve activations on meeting a hydrogen demand, and to what extent this might limit the value of ancillary services for electrolyzers. Further, there is a lack of comparison between specific reserves and how their bid structure affect their value for an electrolyzer.
Based on these identified research gaps, we summarize our main contributions as follows:
\begin{itemize}
    \item A value analysis of participating in the FCR and mFRR markets using historical prices from the years 2021-2023,
    \item An analytical derivation of the opportunity cost (i.e., willingness to participate) of the electrolyzer providing various reserve services,
    \item An investigation of the impact of reserve activation intensity on the hydrogen production schedule.
\end{itemize}

The remainder of the paper is structured as follows: Section \ref{prelm} describes electrolyzer physics and target markets. Section \ref{ModelFormulation} presents the bidding strategy model. Section \ref{analytical_results} derive the analytical results. Section \ref{results} discusses numerical results. Section \ref{conclusion} concludes the paper. Finally, four appendices offer additional information and include a nomenclature.

\vspace{2mm}
\section{Preliminaries} \label{prelm}
This section describes the electrolyzer's operational characteristics and the target markets. This forms the foundation for the proposed mixed-integer linear optimization model, which will be presented later in Section \ref{ModelFormulation}.

\subsection{Electrolyzers} \label{Electrolyzers}
Electrolyzers are electrochemical devices that consume electricity and water to split water molecules into hydrogen and oxygen. This study focuses on alkaline electrolyzers, but can also be adopted to polymer electrolyte membrane (PEM) electrolyzers, which are both low-temperature and mature technologies \citep{Gtz2016}. We leave high-temperature solid oxide technologies for future research.

The efficiency of an electrolyzer is a non-linear function of its power consumption, with peak efficiency typically observed at around 30\% of nominal load. To accurately represent this behavior while avoiding non-linearities, we adopt the model proposed by \cite{manuel}. This approach approximates the hydrogen production curve, which represents the relationship between hydrogen output and electricity consumption, using a \textit{piecewise linear function} that requires the inclusion of binary variables. Alternatively, a linear or conic relaxation of the production curve, as suggested in \cite{enrica}, could be utilized. 

Electrolyzers operate in three distinct \textit{operational states}: online, standby, and off. In the online state, the electrolyzer actively produces hydrogen and consumes electricity, and its power setpoint can be controlled within ramp rate limits. However, the electrolyzer has a minimum loading requirement of 10-40\%, as lower loads could lead to unsafe operation due to mixture of hydrogen and oxygen gas \citep{deGroot2022}. For larger modular systems, lower minimum loading rates can be achieved \citep{lesniak2024advancedschedulingelectrolyzermodules}.  In the standby state, the electrolyzer does not produce hydrogen, but consumes a small amount of electricity to remain ready for rapid transition to the online state in seconds or minutes. In the off state, the system is completely shut down, neither consuming electricity nor producing hydrogen, and requires up to several hours to transition back to the online state. 

These operational states are modeled not only to provide a more accurate representation of the electrolyzer but also to ensure that reserves are not offered when the electrolyzer is in the off state, as its startup time is too slow, or in the capacity range below 10\%, where reserve provision is infeasible.
We assume that the electrolyzer has a ramp rate sufficient to fully participate in all the considered markets. However, our proposed formulation can be adjusted to account for ramp rate limitations over a time frame of up to one hour. Beyond this, ramp rate constraints would need to be incorporated into the day-ahead  market bidding process. For more information on the technological specifications of various electrolyzer technologies, we refer to \cite{Gtz2016}.

\subsection{Bilateral Hydrogen Contracts}
We assume that the produced hydrogen is sold through bilateral hydrogen purchase agreements at a fixed hydrogen price, given the limited availability of pool-based hydrogen markets. Electrolyzers are generally designed to meet a (contract-specified) hydrogen demand. In our formulation, this demand is modeled as a set of tube trailers arriving on-site and departing according to a predetermined daily schedule. Once a tube trailer arrives, it serves as gas storage with a defined capacity that must not be exceeded. If no tube trailer is present, or if all available trailers are full, the electrolyzer cannot produce hydrogen due to the lack of off-take capacity.

Another key aspect of hydrogen purchase agreements is the potential inclusion of a minimum demand requirement over a specified time period, such as daily. While this requirement offers certainty for the off-taker, it presents operational challenges for the producer, who may be forced to produce hydrogen under unprofitable electricity price conditions.

If the electrolyzer is not fully supplied by co-located renewable energy sources, it can procure additional electricity from energy markets. In this study, the electrolyzer participates in the day-ahead electricity market to secure all its electricity needs. Notably, the European Commission’s recent delegated act on renewable fuels of non-biological origin (RFNBO) \citep{EU} outlines operational conditions under which electrolyzers producing green hydrogen may procure power directly from the grid. However, this work focuses on a single, stand-alone electrolyzer and assumes that all required electricity is purchased from the grid, independent of the RFNBO conditions.

\subsection{Target Electricity and Ancillary Service Markets}
\label{target_markets}
As outlined in the Energinet reports \citep{Energinet,energinet2}, electrolyzers are expected to participate in multiple frequency-supporting ancillary service markets, offering an opportunity to diversify revenue streams. A comprehensive review of ancillary service markets in Denmark is available in \cite{peter}.

This study focuses on two ancillary services in the western Denmark bidding zone (DK1), which is part of the synchronous grid of continental Europe: FCR (primary reserve) and mFRR (tertiary reserve). We leave the secondary reserve, the automatic Frequency Restoration Reserve, out of this analysis, as its current structure is similar to that of mFRR, and as its Danish market has recently undergone major restructuring.  These reserve products are also available in other zones within the continental European synchronous grid and are increasingly integrated into common European market platforms for capacity and energy trading \citep{NBM}. In the following, we briefly describe the target markets for electrolyzers in Denmark, presented in chronological order of market closure. We assume the electrolyzer acts as a price-taker in all markets.

\textit{FCR market}:
FCR is the fastest frequency service in the continental European synchronous area and is used mainly for smaller deviations in frequency ($\pm 0.1$ Hz). The transmission system operator procures reserves (in MW) in the FCR capacity market that can be automatically activated in real-time to respond to frequency deviations within seconds. \cite{Tuinema}, \cite{Dozein2021} and \cite{Ribeiro2023} demonstrate that electrolyzers have response times that are potentially faster than conventional generators and are therefore capable of participation in the FCR.
The FCR product is symmetric, meaning the bid capacity is reserved for both upward and downward regulation. In DK1, the FCR capacity market closes at 8:00 AM the day before delivery, with bids submitted for 4-hour blocks. Each bid specifies the reserved capacity (MW) and a price (per MW per hour), and the market operates with uniform pricing.

\textit{mFRR market}: The mFRR is the slowest reserve in the continental European synchronous area and is used upon larger deviations in frequency that cannot be mitigated through the primary or secondary reserves alone. The mFRR market comprises two stages: capacity procurement and energy activation. The mFRR capacity market closes at 9:30 AM the day before delivery.  The reserve is asymmetric, meaning that separate bids can be offered in the up- and down directions, with separate clearing prices in each direction. Bids specify reserved capacity (MW) for 1-hour intervals and a price (per MW per hour),\footnote{\label{note1} Note that the mFRR and day-ahead markets will transition to 15-minute intervals in 2025.} and the market operates with uniform pricing. Once a bid is accepted, the electrolyzer must submit an activation price for the energy activation stage. Activation follows a merit-order curve and is remunerated based on the uniform balancing price. The reserved capacities can be activated within 15 minutes and the activation signal lasts 15 minutes.

\textit{Day-ahead electricity market}:
Electrolyzers can procure electricity for hydrogen production through the day-ahead market. Separate bids of prices and quantities can be placed for each hour. In Denmark, the day-ahead market closes at 12:00 PM the day before delivery and operates with uniform pricing and 1-hour intervals for bids.$^{\ref{note1}}$ These bids must account for reserve obligations from the FCR and mFRR markets.

\textit{Reserve activation stage:}
FCR is automatically activated when a deviation in the nominal frequency of 50 Hz is detected. In contrast, mFRR is manually activated by the transmission system operator based on system conditions. When bidding for reserve capacity a day before delivery, the exact volume of activated capacity remains uncertain. The actual energy activated is determined ex-post and compensated based on the balancing price of the relevant hour. While energy activation of FCR is generally negligible, the activation of mFRR can be significant, as we show later in Section \ref{sec:histrotical_activation}.
\vspace{-2mm}
\section{Model Formulation} 
\label{ModelFormulation}
\begin{newtext}
This section presents a deterministic, mixed-integer linear program aimed at maximizing the total profit of an electrolyzer operating in the day-ahead, FCR, and mFRR markets, while selling hydrogen at a fixed price with a minimum production requirement. Before presenting the formulation, a brief summary of the most relevant assumptions is given. 

The formulation assumes a low-temperature and mature electrolyzer technology, such as alkaline or PEM electrolyzers. It is further assumed that the electrolyzer has ramp rates high enough to qualify its full feasible capacity range for either the FCR or mFRR markets. Moreover, the model assumes perfect price information. While in reality the three considered markets close at different times, the formulation adopts a single gate closure. Finally, the formulation does not account for potential revenues from the balancing stage at the time of bidding. We argue that, unlike the day-ahead price, the imbalance price is very challenging to forecast in the day-ahead stage, and the electrolyzer would not have access to this information in practice.
\end{newtext}

\textit{Notational remark}: Hereafter, lower-case symbols are variables, whereas upper-case and Greek symbols are parameters. In addition, symbol $z$ is used for various binary variables. The upper case $Z$ notes the objective value.

\subsection{Objective Function}
\begin{newtext}
The electrolyzer operates to maximize its own profits in the day-ahead, mFRR and FCR markets following the objective function $Z$:
\end{newtext}

\begin{align}
     \underset{\mathcal{X}}{ \operatorname{max} } \hspace{0.2cm} 
 Z &= \sum_{t \in \mathcal{T}} \Big(
h^{\rm{e}}_{t} \lambda^{\rm{h}}  
+  r^{\rm m \uparrow}_{t}  \lambda^{\rm{m\uparrow}}_{t}  
+   r^{\rm m \downarrow}_{t} \lambda^{\rm{m\downarrow}}_{t} 
- p^{\rm{DA}}_{t} \lambda^{\rm{DA}}_{t}\Big) 
\notag \\ & +  \sum_{i \in \mathcal{I}}   r^{\rm F}_{i} \lambda^{\rm{F}}_{i} \Delta T^{\rm{F}}, 
 \label{D2_obj} 
\end{align}
where $\mathcal{X}$ represents the set of variables, which will be defined later in this section. The index $t \in \mathcal{T}$ denotes 1-hour time intervals, and $i \in \mathcal{I}$ corresponds to 4-hour periods, both necessary for modeling the FCR market. The first term in the objective function represents the revenue from selling hydrogen, which is the product of hydrogen production $h^{\rm{e}}_{t}$ and the fixed hydrogen price $\lambda^{\rm{h}}$. The second and third terms correspond to mFRR capacity bids, split into upward ($r^{\rm m \uparrow}_{t}$) and downward ($r^{\rm m \downarrow}_{t}$) reserves, sold at forecasted prices $\lambda^{\rm{m\uparrow}}_{t}$ and $\lambda^{\rm{m\downarrow}}_{t}$, respectively. The fourth term accounts for the cost of procuring power from the day-ahead market, considering the bid quantity $p^{\rm{DA}}_{t}$ and forecasted price $\lambda^{\rm{DA}}_{t}$. Finally, the last term represents the profit from selling FCR capacity, where the quantity $r^{\rm F}_{i}$ is bid for a 4-hour interval $i$ at forecasted price $\lambda^{\rm{F}}_{i}$, with $\Delta T^{\rm{F}} = 4$.

At the time of the bidding decision, accurately forecasting the hourly balancing prices for the day of delivery and the energy to be activated through FCR or mFRR services is challenging. As a result, during the bidding stage, we neglect the potential revenues or costs associated with real-time activation and focus on maximizing the total profit from the day-ahead and reserve markets only. However, the impact of activation will be examined in the subsequent ex-post analysis.

The constraints can be categorized into three main groups: electrolyzer constraints, bidding constraints, and activation constraints. These categories will be detailed in the following three subsections.
 

\subsection{Electrolyzer Constraints}
The electrolyzer is modeled with a variable efficiency, which depends on the operational point. This is captured through a piecewise formulation of the hydrogen production curve. The three operational states of the electrolyzer, as described in Section \ref{Electrolyzers}, are also taken into account.

\textit{Piece-wise hydrogen production curve}:
Each segment $s \in \mathcal{S}$ on the hydrogen production curve is defined within a specific range of electrolyzer power consumption, denoted by $\underline{P}_s$ and $\overline{P}_s$. The active segment in an hour depends on the scheduled electrolyzer power consumption $p^{\rm{e}}_{t,s}$ (assuming no reserves are activated):
\begin{subequations}
\begin{align}
& \underline{P}_s  z^{\rm{e}}_{t,s} \leq p^{\rm{e}}_{t,s} \leq  \overline{P}_s  z^{\rm{e}}_{t,s} && \forall  t \in \mathcal{T}, s \in \mathcal{S},  \label{cc}
\end{align}
where $z^{\rm{e}}_{t,s}$ is a binary variable for the segment selection. The rationale for using these binary variables is that the concave relaxation is no longer tight due to the presence of reserve revenues in the objective function. Furthermore, we must ensure that at most one segment can be active at any time, i.e.,
\begin{align}
& \sum_{s \in \mathcal{S}} z^{\rm{e}}_{t,s} \leq 1 && \forall  t \in \mathcal{T}. \label{D2_1_seg} 
\end{align}

The hydrogen production $h^{\rm{e}}_{t}$ is defined as a linear function of the electrolyzer power consumption $p^{\rm{e}}_{t,s}$, with a slope $A_s$ and intercept $B_s$ depending on the active segment, i.e.,
\begin{align}
& h^{\rm{e}}_{t} = \sum_{s \in \mathcal{S}} \Big(A_s p^{\rm{e}}_{t,s} + B_s z^{\rm{e}}_{t,s}\Big) \Delta T  && \forall t \in \mathcal{T}, \label{D2_prod_curve}
\end{align}
%
where $ \Delta T = 1 $. To ensure a coherent formulation, we define an auxiliary variable $ p^{\rm{tot}}_{t} $, which represents the total power consumption across all segments at any given hour $t$:
\begin{align}
        & p^{\rm{tot}}_{t} = \sum_{s \in \mathcal{S}} p^{\rm{e}}_{t,s}   && \forall t \in \mathcal{T}. \label{D2_p_tot}
\end{align}

\textit{Operating states}:
The three operating states of the electrolyzer are defined by one binary variable each, i.e., $z^{\rm{on}}_{t}$, $z^{\rm{sb}}_{t}$, and $z^{\rm{off}}_{t}$. The online state $z^{\rm{on}}_{t}$ is defined such that hydrogen production occurs only in the online state:
\begin{align}
        & z^{\rm{on}}_{t} = \sum_{s \in \mathcal{S}} z^{\rm{e}}_{t,s}  && \forall t \in \mathcal{T}. \label{D2_on_state}
\end{align}

It is also possible to define the online state as the sum of all segment selection variables without defining a new binary variable, but it is defined here for clarity. The standby state $z^{\rm{sb}}_{t}$ is defined as:
\begin{align}
        & p^{\rm{sb}}_{t} =  z^{\rm{sb}}_{t}P^{\rm{sb}}  && \forall t \in \mathcal{T}, \label{D2_sb_state}
\end{align}
stating that there is no hydrogen production, but the electrolyzer is still consuming power $p^{\rm{sb}}_{t}$ to keep the system pressurized and heated for a quick start-up of production. The power consumption is assumed constant at $P^{\rm{sb}}$.

In the case of a complete shutdown of the electrolyzer and the surrounding system, represented by $z^{\rm{off}}_{t}$, a minimum downtime is required due to the time it takes to reheat and pressurize the system. The downtime is defined by the $t$-dependent set $\mathcal{N}_t$ as:
\begin{align}
    & z^{\rm{off}}_{t} - z^{\rm{off}}_{t-1} \leq z^{\rm{off}}_{n}
        && \forall t \in \mathcal{T} \setminus \{1\},  n \in \mathcal{N}_t,  \label{D2_down_time}
\end{align}
where $\mathcal{N}_t$ is defined as the $N^{\rm{down}}$ next periods from $t$.
Finally, the electrolyzer must only be in one state at any given hour, i.e.,
\begin{align}
    & z^{\rm{on}}_{t} + z^{\rm{sb}}_{t} + z^{\rm{off}}_{t} = 1  
        && \forall t \in \mathcal{T}.  \label{D2_one_state}
\end{align}
\end{subequations}


\subsection{Bidding Constraints} 
The bids into the day-ahead, FCR, and mFRR capacity markets are all limited by the scheduled electrolyzer consumption and power capacity. Given our assumption that the electrolyzer is a price-taker in all markets, the electrolyzer places a quantity bid based on a forecasted (or known) price. In addition, we assume that there is no information gain between the gate closure times of the various markets and, therefore, that all bids can be placed simultaneously.

\textit{Day-ahead quantity bids}:
The day-ahead bid is defined as the total scheduled consumption of the electrolyzer, including the stand-by consumption:
\begin{subequations}
 \begin{align}
& p^{\rm{DA}}_{t} = p^{\rm{tot}}_{t} + p^{\rm{sb}}_{t} 
        && \forall t \in \mathcal{T}.  \label{D2_DA_bid}
 \end{align}

Note that the power consumption of the compressor and other system components is assumed to be embodied within the hydrogen production curve of the electrolyzer, therefore it is not explicitly indicated in \eqref{D2_DA_bid}. 

\textit{Reserve quantity bids}:
As the FCR market is symmetric, i.e., the identical capacity $r^{\rm F}_{i}$ is reserved for activation in both up and down directions, it must be limited both by the ability to up-regulate and down-regulate the power consumption of the electrolyzer. The mFRR reserve is, however, split between the upward and downward reserve, hence each reserve bid $r^{\rm m \downarrow}_{t}$ and $r^{\rm m \uparrow}_{t}$ is limited in its respective direction only. As the electrolyzer can participate in both markets simultaneously, it is the sum of reserves in each direction that is constrained. For the downward reserve, this is simply limited by the difference between the electrolyzer capacity $C^{\rm{e}}$ and the scheduled power consumption $p^{\rm{tot}}_{t}$. Further, as the electrolyzer has a cold start time of up to 2 hours, the reserve cannot be provided if the electrolyzer is shut off. For the downward reserve, we enforce:
 \begin{align}
& r^{\rm F}_{i} +  r^{\rm m \downarrow}_{t}  \leq C^{\rm{e}} \big(1 - z^{\rm{off}}_{t}\big)  - p^{\rm{tot}}_{t}
        && \forall t \in \mathcal{T}_i,  i \in \mathcal{I},  \label{D2_both_bids_dw}
\end{align}
where $\mathcal{T}_i$ is the set of hours $t$ within the $4$-hour interval $i$. For the upward reserve, the total bids are limited by the difference in the scheduled power consumption and the minimum allowed electrolyzer power $\underline{E}$, i.e.,  
 \begin{align}
& r^{\rm F}_{i} +  r^{\rm m \uparrow}_{t}  \leq   p^{\rm{tot}}_{t} -  \underline{E} \big(1 - z^{\rm{off}}_{t}\big) 
        && \forall t \in \mathcal{T}_i,  i \in \mathcal{I}, \label{D2_both_bids_up}
 \end{align}
where the off-state binary variable $z^{\rm{off}}_{t}$ is included to avoid infeasibility in the case of no consumption.

\textit{Minimum and maximum bid sizes}: Both FCR and mFRR capacity markets operate with minimum and maximum allowed bid sizes. These limits are enforced by the following constraints, with the binary variable $z^{\rm{F}}_{i}$ enabling an FCR bid in the range of $\underline{P}^{\rm{F}}$ and $\overline{P}^{\rm{F}}$. The binary variables $z^{\rm{m\uparrow}}_{t}$ and $z^{\rm{m\downarrow}}_{t}$ serve similarly for the mFRR bids:
\begin{align}
& \underline{P}^{\rm{F}}  z^{\rm{F}}_{i} \leq r^{\rm F}_{i} \leq  \overline{P}^{\rm{F}}  z^{\rm{F}}_{i} && \forall  i \in \mathcal{I} \label{D2_size_bids_FCR}\\
& \underline{P}^{\rm{m}}  z^{\rm{m\downarrow}}_{t} \leq r^{\rm m \downarrow}_{t} \leq  \overline{P}^{\rm{m}}  z^{\rm{m\downarrow}}_{t} && \forall  t \in \mathcal{T}  \label{D2_size_bids_mFRR_dw} \\
& \underline{P}^{\rm{m}}  z^{\rm{m\uparrow}}_{t} \leq r^{\rm m \uparrow}_{t} \leq  \overline{P}^{\rm{m}}  z^{\rm{m\uparrow}}_{t} && \forall  t \in \mathcal{T}.  \label{D2_size_bids_mFRR_up} 
\end{align}
\end{subequations}

If the electrolyzer ramp-rate is too slow to deliver its full capacity within the required time for a reserve, the upper limits of these constraints can be tightened.

\subsection{Activation Constraints}
The capacity reserved in ancillary service markets can be activated upon deviations of the grid frequency. This causes deviations in the hydrogen production, which might be infeasible in a system where there are upper and lower limits to the possible hydrogen off-take. For the FCR bids, it is assumed that the activation is symmetric in both directions, so that the total difference from the scheduled hydrogen production due to activation is negligible. This assumption is tested later in Section \ref{results}. For the mFRR bids, however, the activation could be energy-intensive, and a large deviation in the scheduled production could yield full tube-trailers or unmet minimum demand requirements. Assuming an expected upward mFRR activation $\alpha^{\rm{act \uparrow}}_{t}$ or a downward mFRR activation $\alpha^{\rm{act \downarrow}}_{t}$ at time period $t$, a set of constraints is defined to ensure that the capacity bids are feasible in terms of deviations in hydrogen production. Note that $\alpha^{\rm{act (.)}}_{t} = 1$ refers to the full capacity activation for one hour.

\textit{Upward mFRR activation}:
In this case, the actual power consumption  of the electrolyzer is defined as:
\begin{subequations} \label{activation}
\begin{align}
    & \sum_{s \in \mathcal{S}} p^{\rm{act \uparrow}}_{t,s} = p^{\rm{tot}}_{t} - r^{\rm m \uparrow}_{t}  \alpha^{\rm{act \uparrow}}_{t}   && \forall t \in \mathcal{T}.  \label{D2_P_act_up_def}  
\end{align}

Note that upward activation can only lead to a \textit{decrease} in consumed energy, hence we are only interested in the lower limits on the produced hydrogen. 

Let the binary variable $z^{\rm{act \uparrow}}_{t,s}$ define the range for where a segment on the production curve is active based on the activated power consumption. By this, we restrict $p^{\rm{act \uparrow}}_{t,s}$ by:
\begin{align}
    & \underline{P}_s  z^{\rm{act \uparrow}}_{t,s} \leq p^{\rm{act \uparrow}}_{t,s} \leq  \overline{P}_s  z^{\rm{act \uparrow}}_{t,s} && \forall  t \in \mathcal{T}, s \in \mathcal{S}, \label{D2_P_act_up_seg} 
\end{align}
while only one segment is allowed to be active:
\begin{align}
    & \sum_{s \in \mathcal{S}} z^{\rm{act \uparrow}}_{t,s} \leq 1 && \forall  t \in \mathcal{T}. \label{D2_1_act_up_seg}  
\end{align}

The hydrogen production following an upward activation can then be calculated based on the power consumption and corresponding segment on the production curve, i.e., 
\begin{align}
& h^{\rm{act \uparrow}}_{t}   =  \sum_{s \in \mathcal{S}} \big(A_s p^{\rm{act \uparrow}}_{t,s} + B_s z^{\rm{act \uparrow}}_{t,s}\big) \Delta T      && \forall t \in \mathcal{T}.  \label{D2_H_act_up_def} 
 \end{align}
 
Even if the full activation of mFRR occurs,  the hydrogen production is still constrained by the minimum hydrogen demand $\underline{H}$ over time period $\mathcal{T}$, i.e., 
 \begin{align}
&  \sum_{t \in \mathcal{T}} h^{\rm{act \uparrow}}_{t} \geq \underline{H}
. 
        && \label{D2_H_act_up_daily}
\end{align}

Note that the hydrogen demand is a minimum demand, and can be exceeded while respecting other physical constraints.

\textit{Downward mFRR activation}:
In this case, the actual power consumption  of the electrolyzer is defied as:
\begin{align}
    & \sum_{s \in \mathcal{S}} p^{\rm{act \downarrow}}_{t,s}   =  p^{\rm{tot}}_{t} + r^{\rm m \downarrow}_{t}  \alpha^{\rm{act \downarrow}}_{t}   && \forall t \in \mathcal{T}. \label{D2_P_act_dw_def} 
\end{align}

Note again that downward activation can only lead to an \textit{increase} in consumed energy, hence we are only interested in the upper limits on the produced hydrogen.

We define the binary variable $z^{\rm{act \downarrow}}_{t,s}$ to indicate the corresponding production curve segment. Similar to the upward activation constraints, we enforce:
\begin{align}
    & \underline{P}_s  z^{\rm{act \downarrow}}_{t,s} \leq p^{\rm{act \downarrow}}_{t,s} \leq  \overline{P}_s  z^{\rm{act \downarrow}}_{t,s} && \forall  t \in \mathcal{T}, s \in \mathcal{S} \label{D2_P_act_dw_seg}  \\
    & \sum_{s \in \mathcal{S}} z^{\rm{act \downarrow}}_{t,s} \leq 1 && \forall  t \in \mathcal{T}. \label{D2_1_act_dw_seg}  
\end{align}

The actual hydrogen production following a downward activation can then be calculated as:
\begin{align}
& h^{\rm{act \downarrow}}_{t}   =  \sum_{s \in \mathcal{S}} \Big(A_s p^{\rm{act \downarrow}}_{t,s} + B_s z^{\rm{act \downarrow}}_{t,s}\Big) \Delta T      && \forall t \in \mathcal{T}.  \label{D2_H_act_dw_def} 
 \end{align}
 
In this case, the additional hydrogen production should respect the upper limit of available capacity of tube-trailers on site. Let $d \in \mathcal{D}$ a set of tube-trailers. A set of auxiliary variables is defined tracking the dispensed hydrogen  $h^{\rm{d}}_{t,d}$ and state of tube-trailers for the activated hydrogen flow $s^{\rm{d}}_{t,d}$. These flows are then constrained by the dispenser capacity $C^{\rm{d}}$, the tube-trailer availability indicated by the binary parameter $\zeta_{t,d}$,  and their capacity $S^{\rm{d}}$: 
\begin{align}
& h^{\rm{act \downarrow}}_{t} = \sum_{d \in \mathcal{D}}     h^{\rm{d}}_{t,d} 
        && \forall t \in \mathcal{T} \label{D2_H_act_dw_disp} \\
& h^{\rm{d}}_{t,d} \leq C^{\rm{d}} \zeta_{t,d} 
        && \forall t \in \mathcal{T}, d \in \mathcal{D} \label{D2_H_act_dw_tt} \\
& s^{\rm{d}}_{t,d} = h^{\rm{d}}_{t, d}  
        && \forall t = 1,  d \in \mathcal{D} \label{D2_H_act_dw_s0} \\
& s^{\rm{d}}_{t,d} = s^{\rm{d}}_{t-1,d} + h^{\rm{d}}_{t, d} 
        && \forall t \in \mathcal{T} \setminus \{1\},  d \in \mathcal{D}  \label{D2_H_act_dw_st} \\
& s^{\rm{d}}_{t,d} \leq S^{\rm{d}}
        && \forall t \in \mathcal{T}, d \in \mathcal{D}. \label{D2_H_act_dw_smax}
\end{align}
\end{subequations}

With separate sets of constraints for upward and downward activated hydrogen flows, it is ensured that all market bids are feasible for any assumed expected activation. This is intrinsically a conservative approach, as activation in one direction can be fully or partially canceled out by activation in the other, which is not accounted for in this formulation. However, accounting for this effect requires information regarding the timing and frequency of activation, as the number of consecutive hours of activation in one direction is of significance when constrained by tube trailer capacity and minimum demand requirements. An alternative approach is to only apply the activation constraints on the \textit{net} activation. In the proposed formulation, the risk can be adjusted through the selection of expected activation, $\alpha^{\rm{act \uparrow}}_{t}$ and $\alpha^{\rm{act \downarrow}}_{t}$. If set to zero, the model does not take activation into account at all. If set to one, the bids are robust in terms of activation in any direction. 
\begin{newtext}
    The electrolyzer operator can select the $\alpha$-parameter based on their level of risk aversion. For example, if failing to meet the hydrogen demand results in a high penalty, choosing $\alpha = 1$ may be appropriate.
\end{newtext}


As a summary, the proposed mixed-integer linear optimization problem is formulated in  \eqref{D2_obj}-\eqref{activation}, with the set of variables
    $\mathcal{X}$ =  $\{ h^{\rm{act \uparrow}}_{t}$, $h^{\rm{act \downarrow}}_{t}$, $h^{\rm{d}}_{t,d}$, $h^{\rm{e}}_{t}$,  $p^{\rm{act \uparrow}}_{t,s}$, $p^{\rm{act \downarrow}}_{t,s}$, $p^{\rm{DA}}_{t}$, $p^{\rm{e}}_{t,s}$, 
    $p^{\rm{sb}}_{t}$,  $p^{\rm{tot}}_{t}$, 
    $r^{\rm F}_{i}$, $r^{\rm m \uparrow}_{t}$, $r^{\rm m \downarrow}_{t}$, $s^{\rm{d}}_{t,d}$, $z^{\rm{act \uparrow}}_{t,s}$, $z^{\rm{act \downarrow}}_{t,s}$, $z^{\rm{e}}_{t,s}$, $z^{\rm{F}}_{i}$, $z^{\rm{m \uparrow}}_{t}$, $z^{\rm{m \downarrow}}_{t}$, $z^{\rm{off}}_{t}$, $z^{\rm{on}}_{t}$, $z^{\rm{sb}}_{t}\}$. Note that   $\boldsymbol{z}$ is the set of binary variables, and the rest  are non-negative continuous variables.

\section{Analytical Results: When Is It Optimal to Offer a Reserve Capacity?}
\label{analytical_results}
In this section, we analytically derive the reserve price that the electrolyzer must charge to offer reserve capacity, based on the reserve bidding structure. We consider three different reserve structures:
\begin{itemize}
    \item \textbf{mFRR:} Hourly bid, asymmetric reserve. This structure is equivalent to the current mFRR structure in DK1.
    \item \textbf{Relaxed FCR:} Hourly bid, symmetric reserve. This is similar to the current FCR structure in DK1, but with the four-hour bid duration relaxed.
    \item \textbf{FCR:} Four-hour bid, symmetric reserve. This is equivalent to the current FCR structure in DK1.
\end{itemize}

The bid prices for each reserve structure are derived from the \textit{opportunity cost} of providing reserve capacity. We show that this cost depends on the day-ahead price and the profit the electrolyzer achieves when optimizing its schedule based solely on the day-ahead price. The derived bid curves, representing the minimum reserve price as a function of the reserved quantity, take the form of stepwise increasing curves. We illustrate these curves using a stylized example.


\subsection{Derivation of Reserve Bid Curves}  
Let $\tilde{p}^{\rm DA*}_{t}$ denote the electrolyzer's optimal day-ahead schedule at time $t$ when it does not participate in any reserve markets. This optimal quantity is determined by the optimization model \eqref{D2_obj}-\eqref{activation} given the forecast price for the day-ahead market. For simplicity, we assume that the electrolyzer remains online, ignoring standby consumption and start-up time. Consequently, the electrolyzer setpoint $\tilde{p}^{\rm tot *}_{t}$ is equal to $\tilde{p}^{\rm DA*}_{t}$.

If the electrolyzer intends to reserve a capacity in the upward or downward direction for one hour, it might have to change its setpoint from $\tilde{p}^{\rm DA*}_{t}$, depending on the magnitude of the reserved capacity. Let ${\rm p}^{\rm \downarrow}(r_{t})$ denote the electrolyzer's new setpoint as a function of offering a downward reserve capacity $r_{t}$:
\begin{subequations}
\label{eq:ely_setpoint_with_reserve_both}
\begin{align}
\label{eq:down_resulting_setpoint}
    & {\rm p}^{\rm \downarrow}(r_{t}) = 
    \begin{cases}
        \tilde{p}^{\rm DA*}_{t} & \text{if} \quad r_{t} \leq \bar{P}- \tilde{p}^{\rm DA*}_{t}, \\ 
        \bar{P}- r_{t} & \text{if else}.
    \end{cases}
\end{align}

Similarly, let ${\rm p}^{\rm  \uparrow}(r_{t})$ be the electrolyzer's new setpoint as a function of offering an upward reserve capacity $r_{t}$:
\begin{align}
\label{eq:up_resulting_setpoint}
    & {\rm p}^{\rm  \uparrow}(r_{t}) = 
    \begin{cases}
        \tilde{p}^{\rm DA*}_{t} & \text{if} \quad r_{t} \leq  \tilde{p}^{\rm DA*}_{t} - \underbar{P}, \\ 
        \underbar{P} + r_{t} & \text{if else}.
    \end{cases}
\end{align}
\end{subequations}

Equations \eqref{eq:ely_setpoint_with_reserve_both} imply that the electrolyzer can bid any capacity above  $\tilde{p}^{\rm DA*}_{t}$ to downward reserves, and any capacity below to upward reserves, without diverting from this setpoint. 

We now derive the profit that the electrolyzer generates by operating at a given setpoint, excluding any revenue from reserved capacities. The electrolyzer's profit at time $t$, denoted by $\tilde{\Gamma}_{t}$, is determined from the optimal day-ahead bid $\tilde{p}^{\rm DA*}_{t}$, assuming the electrolyzer is not participating in any reserve markets. This profit can be calculated from the objective function of the optimization model \eqref{D2_obj}-\eqref{activation}:
\begin{align}
    \tilde{\Gamma}_{t} = \lambda^{\rm h}  \tilde{h}^{\rm e*}_{t} - \lambda^{\rm DA}_{t}  \tilde{p}^{\rm DA*}_{t},
\end{align}
where $\tilde{h}^{\rm e*}_{t}$ represents the (approximated) hydrogen produced when operating at the $\tilde{p}^{\rm DA*}_{t}$ setpoint. This value can either be directly taken from the optimization model \eqref{D2_obj}-\eqref{activation} or computed as follows:
\begin{align}
\label{eq:hydrogen_production_approx}
    \tilde{h}^{\rm e*}_{t} = {\rm h}(\tilde{p}^{\rm DA*}_{t} ) = \sum_{s \in \mathcal{S}} (A_s \tilde{p}^{\rm DA*}_{t} + B_s)\mathbf{1}_{\left(\tilde{p}^{\rm DA*}_{t} \in [\underline{P}_s, \overline{P}_s ]\right)}, 
\end{align}
where $\mathbf{1}_{(.)}$ is the indicator function, which takes a value of 1 when $\tilde{p}^{\rm DA*}_{t} \in [\underline{P}_s, \overline{P}_s ]$, and 0 otherwise. The function $\rm{h}(\cdot)$ can thus be used to compute the hydrogen production corresponding to any electrolyzer setpoint.

Using \eqref{eq:hydrogen_production_approx}, we can define the profit (excluding reserve market revenues) as a function of any electrolyzer setpoint obtained from \eqref{eq:ely_setpoint_with_reserve_both}. The profit, $\tilde{\Gamma}(p^{\rm (\cdot)}_{t})$, corresponding to the electrolyzer's revised setpoint $p^{\rm (\cdot)}_{t}$ due to the provision of reserve capacity, is given by:
\begin{align}
    \tilde{\rm \Gamma}{(p^{\rm (\cdot)}_{t})} = \lambda^{\rm h}   {\rm h}(p^{\rm (\cdot)}_{t} ) - \lambda^{\rm DA}_{t}  p^{\rm (\cdot)}_{t}  .
\end{align}

If the profit $\tilde{\Gamma}(p^{\rm (\cdot)}_{t})$ from providing a reserve is lower than the optimal day-ahead profit, an opportunity cost is associated with reserving the capacity. In such a case, the electrolyzer will only reserve capacity if the revenue from the reserve meets or exceeds this opportunity cost. The opportunity cost ${\rm C}^{(\cdot)}(r_{t})$ incurred by the electrolyzer when providing a reserve capacity can be determined as:
\begin{subequations}
\label{eq:OC_of_reserve_both}
\begin{align}
\label{eq:OC_down}
    {\rm C}^{\downarrow}(r_{t}) &= \tilde{\Gamma}(\tilde{p}^{\rm DA*}_{t}) - \tilde{\Gamma}({\rm p}^{\rm \downarrow}(r_{t})), \\
    \label{eq:OC_up}
    {\rm C}^{\uparrow}(r_{t}) &= \tilde{\Gamma}(\tilde{p}^{\rm DA*}_{t}) - \tilde{\Gamma}({\rm p}^{\rm  \uparrow}(r_{t})),
\end{align}
\end{subequations}
where the opportunity cost is the difference between the profit achieved with the optimal day-ahead setpoint, $\tilde{p}^{\rm DA*}_{t}$, and the profit achieved with the setpoint resulting from offering reserve capacity, either ${\rm p}^{\rm \downarrow}(r_{t})$ or ${\rm p}^{\rm  \uparrow}(r_{t})$. This derivation of the reserve-induced opportunity cost relies on two assumptions:
\begin{itemize}
    \item The reserve provided will not be activated, and
    \item For downward reserves, the minimum demand constraint, as enforced in  (\ref{D2_H_act_up_daily}), must not be binding.
\end{itemize}

Since the hydrogen production curve in \eqref{eq:hydrogen_production_approx} is approximated as piecewise linear, the opportunity cost curves in \eqref{eq:OC_of_reserve_both} are also piecewise linear. According to the definition of ${\rm p}^{\rm \downarrow}$ in \eqref{eq:down_resulting_setpoint}, \eqref{eq:OC_down} evaluates to zero for any reserved capacity up to $\bar{P} - \tilde{p}^{\rm DA*}_{t}$. Therefore, depending on the setpoint $\tilde{p}^{\rm DA*}_{t}$, the electrolyzer may be willing to offer downward reserve capacity for any positive reserve price. Similarly, in the upward direction, the electrolyzer is willing to offer reserve capacity up to $\tilde{p}^{\rm DA*}_{t} - \underline{P}$ for any positive reserve price.

To offer capacity beyond this, the electrolyzer requires a minimum price to compensate for the reduced profit resulting from changing its day-ahead setpoint. This minimum price, i.e., the electrolyzer's bid price, noted $\pi$, is determined by taking the derivative of the opportunity cost functions in \eqref{eq:OC_of_reserve_both}:
\begin{subequations}
\label{eq:price_of_bidcurve_reserve_both}
\begin{align}
\label{eq:min_price_down}
    \pi^{\downarrow}(r_{t}) = {\rm C}^{\downarrow \prime }  (r_{t}), \\
    \label{eq:min_price_up}
\pi^{\uparrow}(r_{t}) = {\rm C}^{\uparrow \prime } (r_{t}).
\end{align}
\end{subequations}

The mFRR, as explained in Section \ref{target_markets}, is provided in hourly blocks and separated into the upward and downward directions. Therefore, to provide an mFRR service, the electrolyzer requires a minimum price from the reserve, as given by (\ref{eq:min_price_down}) and (\ref{eq:min_price_up}), for the downward and upward directions, respectively. Equations (\ref{eq:price_of_bidcurve_reserve_both}) can thus be used to form the electrolyzer's bid curve for mFRR provision, assuming perfect information of the day-ahead price.

The FCR, however, is a symmetric reserve, meaning the reserve provided might be activated in either direction. The bid price for providing a reserve in both directions, $\pi^{\rm \downarrow \uparrow}$, is given by the maximum of the upward and downward bid prices of (\ref{eq:price_of_bidcurve_reserve_both}):
\begin{align}
\label{eq:min_price_symmetric}
    \pi^{\rm \downarrow \uparrow}(r_{t}) = \text{max}\left\{\pi^{\downarrow}(r_{t}), \pi^{\uparrow}(r_{t})\right\}.
\end{align}

Additionally, the derivation so far assumes a reserve block duration of one hour, while the current FCR block duration in DK1 is four hours. When providing reserve capacity for four hours, the minimum revenue must account for the difference in revenue achieved through the day-ahead market over all four hours. Therefore, the bid price for providing a symmetric reserve capacity over a time block $i$ is determined as the average price of all hours within that block:
\begin{align}
\label{eq:min_price_4h}
    \pi^{\rm \downarrow \uparrow, 4h}(r_i) & = \frac{1}{|\mathcal{T}_i |}\sum_{t \in \mathcal{T}_i  }\pi^{\rm \downarrow \uparrow}_{t}(r_i),
\end{align}
where $\mathcal{T}_i$ is the set of hours belonging to block $i$, which in this case has a length of four. Equation (\ref{eq:min_price_4h}) forms the electrolyzer's bid curve for FCR provision assuming perfect information of the day-ahead price, while  (\ref{eq:min_price_symmetric}) forms the FCR bid curve under the same assumption for a relaxed one-hour block structure.

\subsection{Stylized Example}
We illustrate the derived reserve bid curves with a constructed example that spans eight hours. In this example, all reserves have a constant price of 11 EUR/MW/h, while the day-ahead price changes per hour, ranging from 3 to 106 EUR/MWh.

\begin{figure}[t]
    \centering
    \includegraphics[width=\linewidth]{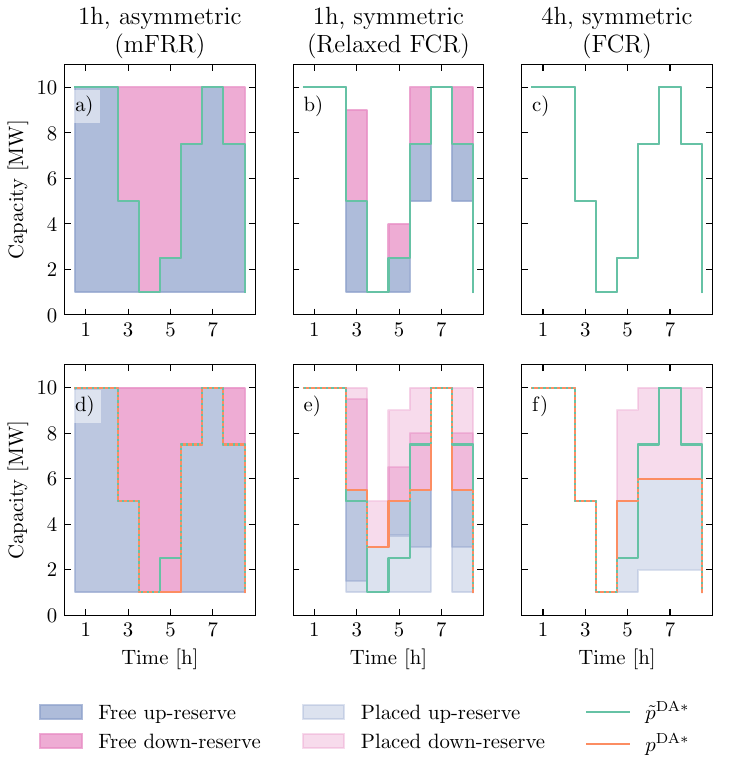}
    \caption{Illustrative example. \textit{Top row:} The amount of reserve capacity the electrolyzer is willing to offer for free, i.e., at any positive reserve price, under different reserve market structures, shown relative to the electrolyzer's day-ahead schedule. \textit{Second row:} The amount of reserves actually placed, considering the actual reserve price, shown relative to the (updated) day-ahead schedule. 
    }
    \label{fig:willingness_free_reserves}
\end{figure}

Figure \ref{fig:willingness_free_reserves} illustrates the reserve capacities the electrolyzer is willing to provide at any positive price (top row) and the actual capacities it places (bottom row) in the stylized example with a constant reserve price. The figure is divided into three columns: the first represents the 1-hour asymmetric reserve structure (mFRR), the second corresponds to the 1-hour symmetric reserve structure (relaxed FCR), and the third depicts the 4-hour block symmetric reserve structure (current FCR). The green line shows the electrolyzer’s optimal schedule based solely on the day-ahead price, while the orange line reflects the updated schedule when optimized for both the day-ahead and reserve prices.

Let us first examine the top row plots (a)-(c). Plot (a) illustrates that the electrolyzer can offer its entire capacity for either mFRR up or down reserve at any positive reserve price. In plot (b), when the reserve structure is symmetric in both directions, the available capacity with no associated opportunity cost decreases. Further constraining the reserve structure to a 4-hour block (plot (c)) leaves no capacity available for free, as the optimal day-ahead price schedule reaches the electrolyzer’s maximum or minimum load in both 4-hour blocks. 


If the reserve prices are sufficiently high, the electrolyzer may adjust its day-ahead price optimal schedule to accommodate larger reserved capacities. This is shown in the second row of Figure \ref{fig:willingness_free_reserves}. In the mFRR structure (plot (d)), only a minor adjustment to the day-ahead schedule occurs in hour 5, where the electrolyzer increases its setpoint to meet the minimum reserve bid size of 2 MW. In the symmetric structure (plot (e)), the amount of capacity placed is reduced compared to plot (d). Notably, no reserve is placed in hours 1, 2, and 7 where the schedule remains at the day-ahead optimized setpoint of maximum capacity. In other hours, the electrolyzer adjusts its schedule to a mid-load position to accommodate more symmetric reserve capacity. In the 4-hour block structure (plot (f)), no reserve bids are placed during the first block of hours 1-4. However, in the second block (hours 5-8), the electrolyzer places more reserve bids than under the 1-hour symmetric structure. This increase is driven by hour 7, where the 4-hour block structure forces the electrolyzer to maintain reserves, even when it incurs a loss in hour 7. However, the total additional revenue during the four hours results in a net positive profit for the electrolyzer. While the electrolyzer schedules more reserves in this block compared to the relaxed FCR structure, its overall profit decreases.

\begin{figure}[t]
    \centering
    \includegraphics[width=\linewidth]{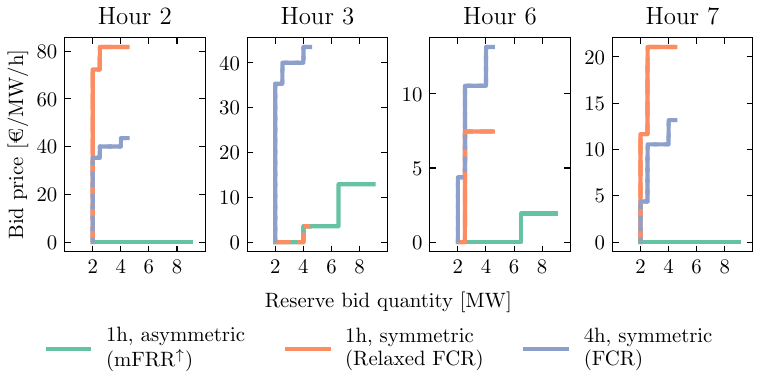}
    \caption{Illustrative example. Reserve bid prices as a function of the reserve capacity under various reserve structures.}
    \label{fig:willingness_prices}
\end{figure}

In this example, the reserve price is constant across all reserves, which means that the differences in the electrolyzer’s willingness to participate in reserves arise solely from changes in the day-ahead price and the differences in reserve structures. To better understand the scheduling results, we examine the reserve bid curves for four specific hours in Figure \ref{fig:willingness_prices}.

In hour 2, the mFRR up reserve price is zero for any reserve quantity because the day-ahead optimal setpoint is at the maximum electrolyzer capacity. Both symmetric reserve structures, regardless of block length, show high bid prices in this hour, as hydrogen production is highly profitable. Reducing the electrolyzer’s setpoint to accommodate a symmetric reserve would therefore incur a high opportunity cost. As the reserve price is lower than the symmetric bid prices, no symmetric reserve is scheduled in hour 2.

In hour 3, the day-ahead optimal setpoint has dropped to a mid-range value. The symmetric bid price under the 1-hour structure is now low enough for the electrolyzer to offer its full capacity for reserve. However, the symmetric bid price for the 4-hour block is still too high because it averages the prices of hours 1-4. As a result, the electrolyzer only schedules a symmetric reserve under the relaxed FCR structure.

In hour 6, the bid price for all reserves is sufficiently low, and the electrolyzer places reserve capacity under all structures. The capacity placed under the 4-hour FCR structure is slightly lower than under the relaxed FCR structure, as the last part of its bid curve is above the example reserve price.

In hour 7, the 4-hour symmetric reserve price is low enough to allow the electrolyzer to offer reserve capacity, while the 1-hour block price remains too high. As mentioned earlier, the 4-hour block price is an average of all four hours. Although the reserve for hour 7 incurs a loss, it is ultimately more profitable to provide the reserve because of the net profit from the other hours in the block.


This example illustrates that the electrolyzer's cost to participate in the reserve markets and consequently its reserve bid price can vary greatly depending on the day-ahead price and the reserve bid structure.

\section{Numerical Results: Case Study of Western Denmark}
\label{results}
 
We apply the proposed optimization model \eqref{D2_obj}-\eqref{activation} to a case study of the western market zone in Denmark, DK1. Section \ref{input-data-section} provides the input data for the case study. Section \ref{activation-unaware} presents the day-ahead scheduling decisions, along with the resulting ex-post profit and hydrogen production following reserve activation, when bidding in an ``activation-unaware" manner (i.e., assuming no activation of reserved capacities, $\alpha = 0$). Finally, we present the day-ahead scheduling and ex-post results when bidding in an ``activation-aware" manner ($\alpha > 0$) in Section \ref{activation-aware}. All source codes are available online \citep{github_AS}.

\subsection{Input Data}
\label{input-data-section}
We present a case study based on historical energy and reserve prices in DK1. This zone is connected to the synchronous grid of continental Europe, unlike the eastern zone, DK2, which is connected to the synchronous Nordic grid. As a result, DK1 shares reserve products with neighboring countries, such as Germany. In addition, the Danish government plans to construct a hydrogen pipeline connecting western Denmark with industrial demands in Germany, making DK1 a likely zone for electrolyzer investments in Denmark \citep{energinet2}.  We focus on the years 2021-2023 to analyze the profits of the electrolyzer. However, due to the availability of frequency data, the majority of our analysis centers on 2023. Since the Danish transmission system operator, Energinet, only procured mFRR \textit{up} quantities in the considered years, we will henceforth refer to ``mFRR up" simply as mFRR. The electrolyzer system input parameters, such as electrolyzer capacity and hydrogen demand, are provided in Table \ref{tab:parameters}. In addition, we assume a daily tube-trailer schedule where all trailers are exchanged at midnight.

\begin{table}[t]
\centering
\caption{Electrolyzer system parameters}
\label{tab:parameters}
\begin{tabular}{llrl}
 \toprule
Electrolyzer capacity & $C^{\rm e}$ & 10     & MW      \\
Minimum load          & $\underline{E}$  & 1      & MW      \\
Standby power         & $P^{\rm{sb}}$      & 0.25   & MW      \\
Base minimum demand           & $\underline{H}$    & 2000   & kg/day  \\
Peak efficiency       & $\bar{\eta}$       &    19   & kg/MW \\
Tube trailer capacity & $C^{\rm{d}}$       & 1000$\times$5 & kg/day   \\
Fixed hydrogen price & $\lambda^{\rm{h}}$       & 5 & EUR/kg   \\\bottomrule  
\end{tabular}
\end{table}

\subsubsection{Historical Prices}
\begin{figure}[t]
    \centering
    \includegraphics[width=\linewidth]{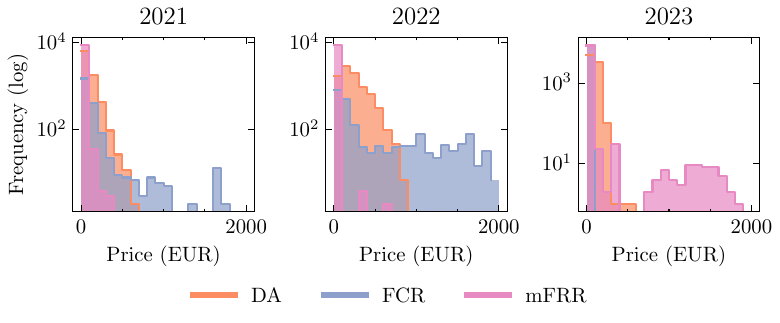}
    \caption{Price distributions for the day-ahead (DA) electricity market, as well as the FCR and mFRR capacity markets in DK1, for the years 2021-2023. Note that the mFRR prices reflect only the up-regulation prices.}
    \label{fig:prices}
\end{figure}

We provide an ex-post analysis \textit{in hindsight}, assuming that the electrolyzer operates with perfect price forecasts. As mentioned earlier, we assume the electrolyzer is a price-taker and does not influence historical prices. Market prices from the day-ahead, balancing, mFRR, and FCR markets can be accessed through the Energi Data Service portal provided by Energinet \citep{energidataservice}. Figure \ref{fig:prices} displays the price distributions for the mFRR, FCR, and day-ahead markets from 2021 to 2023. Note that the y-axis is logarithmic. The price distributions show considerable variation across the years. While 2021 and particularly 2022 saw very high FCR prices, the FCR price drops significantly in 2023. In contrast, mFRR prices are notably higher in 2023 compared to the previous two years. The day-ahead prices were highest in 2022 among the three years analyzed.

\subsubsection{Historical Reserve Activation}
\label{sec:histrotical_activation}
We derive the historical activation of the FCR and mFRR services for the year 2023. Recall that FCR is an automatically activated reserve, where reserve providers must adjust their setpoint based on the frequency measurement at their location. Frequency measurements with a resolution of one second for the year 2023 in the synchronous grid of continental Europe were obtained from the online platform ``Netztransparenz" \citep{netztransparenz}, provided by the four German transmission system operators. The derivation of the FCR activation signal from this frequency is presented in \ref{app:FCR}.

Unlike FCR, the mFRR service is a manually activated reserve, and its activation signal cannot be derived directly from the system frequency. As explained in Section \ref{target_markets}, after the mFRR capacity bids are accepted, mFRR capacity providers must submit an activation price. Upon a system imbalance, the transmission system operator activates mFRR capacities in the balancing stage according to the merit order principle, i.e., starting with the lowest activation price bids. Hence, it can be assumed that if the electrolyzer places an activation price bid lower than or equal to the realized balancing price, $\lambda^{\rm{B}}_{t}$, the reserved capacity would be activated in the event of a system imbalance. Assuming the electrolyzer acts as a price-taker, it would place an activation price bid for upward regulation (i.e., reducing consumption) equal to its opportunity cost, which corresponds to the lost revenue from selling hydrogen. Unlike model \eqref{D2_obj}-\eqref{activation} for optimal bidding purposes, we assume a constant electrolyzer efficiency $\Bar{\eta}$ during the activation stage. Therefore, the opportunity cost, and thus the activation price bid, can be calculated as $\Bar{\eta} \lambda^{\rm{h}}$. 

Although the real activation of mFRR services occurs in 15-minute intervals in Denmark, the relevant historical data are only available at an hourly resolution \citep{energidataservice}. Consequently, we assume that the reserved capacity is fully activated for an entire hour if the described conditions are fulfilled. The full derivation of the mFRR activation signal is presented in \ref{app:mFRR}.

\begin{figure}[t!]
\centering
    \begin{subfigure}[t]{0.49\linewidth}
\includegraphics[height= 1\textwidth]{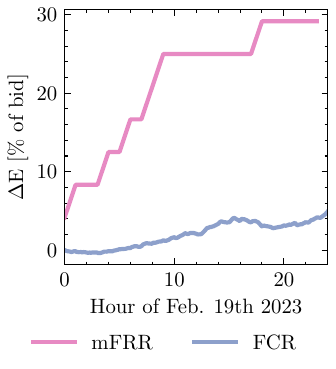}
        \caption{Cumulative activation over a day}
        \label{fig:activation_ex}
    \end{subfigure}%
    \begin{subfigure}[t]{0.49\linewidth}
\includegraphics[height= 1\textwidth]{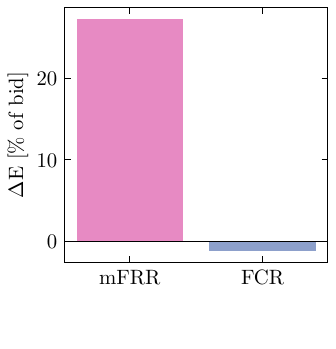}
        \caption{Total activation over a year}
        \label{fig:total_activation}
    \end{subfigure}
    \caption{Historical activation of FCR and mFRR up services over an example day (left) and for the entire year 2023 (right), expressed as a percentage of the reserve bid placed.}
    \label{fig:reserve_activation_23}
\end{figure}

Recall that we hypothesized the activation energy of the mFRR reserve to be significant, whereas that of the FCR reserve is negligible due to its short activation periods and symmetric nature. To evaluate this hypothesis, we analyze the historical activation patterns of each reserve. Figure \ref{fig:reserve_activation_23} presents the \textit{relative} activation of mFRR and FCR reserves, both cumulatively over an example day (Figure \ref{fig:activation_ex}) and in total for the year 2023 (Figure \ref{fig:total_activation}). Here, relative activation, denoted $\Delta E$, refers to the proportion of the reserve bid size that the electrolyzer would have been exposed to during activation, expressed as a percentage. The FCR exhibits activation in both directions, including periods of negative $\Delta E$, due to its symmetric nature. In contrast, the mFRR reserve is divided into separate up and down reserves, with only the up reserve contributing to increasing $\Delta E$. Over the course of 2023, the mFRR up reserve yielded approximately 25 times the net upward activation of the FCR. This validates the assumption that mFRR activation has a significantly larger potential to impact the hydrogen schedule.

\subsubsection{Ex-post Results Upon Reserve Activation}
\label{ex-post-profit-section}
The ex-post profit of the electrolyzer is calculated based on the costs and revenues from the day-ahead electricity, FCR, and mFRR capacity markets, as determined by the objective function \eqref{D2_obj}, as well as the realized activation of reserved capacities, subsequent imbalance payments, and changes in hydrogen production. The ex-post profit, $\Gamma$, is defined as follows:
\begin{align}
    \label{eq:profit_ex_post}
    \rm{\Gamma}  = Z^{*} &- \sum_{t \in \mathcal{T}} \Big[ {\overbrace{\lambda ^{\rm{B}}_{t} ( \Delta E^{\rm{F}}_{t} + \Delta E^{\rm{m \uparrow }}_{t} ) }^{\mathclap{\text{Change in energy consumption}}}} 
    + {\underbrace{\lambda^{\rm{h}} \Delta h_{t}}_{\mathclap{\text{Change in hydrogen production}}}}   \Big],
\end{align}
where the first term, $Z^{*}$, represents the optimal value obtained for the objective function \eqref{D2_obj}. The second term accounts for the additional cost or revenue from the balancing market due to changes in energy consumption resulting from the activation of mFRR, FCR, or both services. If the electrolyzer reduces its energy consumption in response to an activation signal ($\Delta E_{t}^{(.)} < 0$), it will be remunerated at the balancing price $\lambda_t^{\rm{B}}$, and vice versa. The third term reflects the change in revenue from hydrogen production adjustment due to reserve activation, assuming the hydrogen price is fixed at $\lambda^{\rm{h}}$. The change in hydrogen production at hour $t$ is computed as follows:
\begin{align}
    \Delta h_t ={\rm H}\Big(p^{\rm tot}_t  + \Delta P^{\rm{F}}_{t} + \Delta P^{\rm{m \uparrow }}_{t}\Big) - h^{\rm e *}_t,   \ \ \ \forall  t \in \mathcal{T},
\end{align}
where ${\rm H(\cdot)}$ represents the \textit{empirical} (not linearized) hydrogen production as a function of the electrolyzer's power consumption \citep{manuel}, the term $p^{\rm tot}_t + \Delta P^{\rm{F}}_{t} + \Delta P^{\rm{m \uparrow }}_{t}$ denotes the \textit{realized} power consumption of the electrolyzer upon reserve activation, and $h^{\rm e *}_t$ is the optimal hydrogen production schedule at hour $t$, obtained from the optimization model \eqref{D2_obj}-\eqref{activation}. Figure \ref{ex-post} illustrates how we compute the ex-post profit.

\begin{figure} [t]
\centering
\begin{tikzpicture}[node distance = 2.5 cm] 
    \tikzstyle{every node}=[font=\small]

    \tikzstyle{normal} = [rectangle, rounded corners, minimum width=2.5cm, minimum height=1cm,text centered, draw=teal, ultra thick]
    \tikzstyle{long} = [rectangle, rounded corners, minimum width=1.5cm, minimum height=2.5cm,text centered, draw=teal, ultra thick]
    \tikzstyle{wide} = [rectangle, rounded corners, minimum width=3.3cm, minimum height=1cm,text centered, draw=teal, ultra thick]

    \tikzstyle{arrow->} = [->, ultra thick, draw = teal]
    \tikzstyle{arrow<-} = [<-, ultra thick, draw = teal]

    \node [wide] at (0,0) (MILP) {Model \eqref{D2_obj}-\eqref{activation}};
    \node [wide] at (0,-1.5) (input1) { $\lambda_t^{\rm{\{DA,F,m\uparrow,m\downarrow}\}}, \alpha_t^{\rm{act(\cdot)}}$};
    
    \node [long, align=center] at (2.75,-0.75) (output1) {$p_t^{{\rm{DA}}^{*}}$ \\
    $r_t^{{\rm F}^{*}}$ \\
    $r_t^{{\rm m \uparrow}^{*}}$ \\
    $r_t^{{\rm m \downarrow}^{*}}$};
    
    \node [wide, fill = teal!10] at (5.5,0) (expost) {Ex-post profit};
    \node [wide] at (5.5,-1.5) (input3) {$\alpha_t^{\rm{real(\cdot)}}, \lambda_t^{\rm{B}}$};

    \draw [arrow->] (input1) -- (MILP);
    \draw [arrow->]  (MILP) -- (MILP-|output1.west);
    \draw [arrow<-]  (expost) -- (expost-|output1.east);
    \draw [arrow->]  (input3) -- (expost);

\end{tikzpicture}
\vspace{2mm}
    \caption{Overview of the ex-post analysis. Superscript ${}^{*}$ denotes the optimal day-ahead schedules. }\vspace{-2mm}
\label{ex-post}
\end{figure}

\subsection{Activation-Unaware Reserve Participation ($\alpha = 0$)}
\label{activation-unaware}

The primary motivation for an electrolyzer to participate in ancillary service markets is the additional monetary value it can generate through these secondary revenue streams. Therefore, it is crucial that reserve provision increases the electrolyzer's profit in order to make participation in reserve markets worthwhile. To evaluate this, we assess the \textit{maximum} possible increase in profit (upper bound for the profit) that an electrolyzer can achieve through optimal bidding in the FCR and mFRR capacity markets. To determine such a maximum profit, we begin by presenting results for the ex-post profit $\rm{\Gamma}$, calculated as in \eqref{eq:profit_ex_post}. This assumes $\alpha^{\rm act(\cdot)} = 0$, which excludes any potential mFRR activation at the bidding stage.

\begin{figure}[t]
    \centering
    \includegraphics[width=\linewidth]{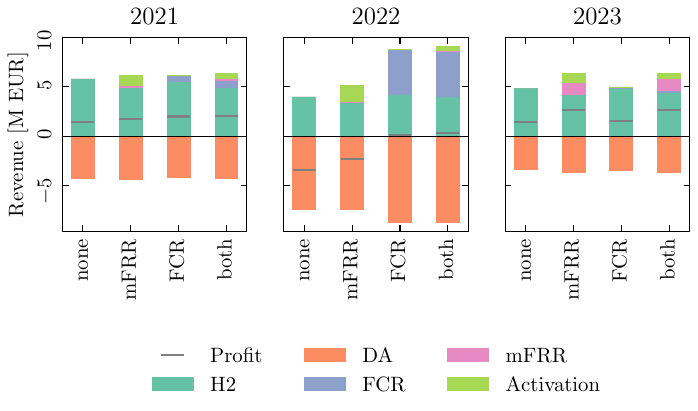}
    \caption{The ex-post profit of the electrolyzer and revenue-cost breakdown, disregarding potential activations in the bidding stage ($\alpha = 0$). DA stands for day-ahead.}
    \label{fig:revenue_streams_all}
\end{figure}

Figure \ref{fig:revenue_streams_all} shows the resulting ex-post profit and the breakdown of revenue-cost streams for the electrolyzer in years 2021-2023. Note that while the ex-post profit is not penalized for failing to meet the minimum hydrogen requirement, it still accounts for changes in hydrogen production and energy consumption, as outlined in \eqref{eq:profit_ex_post}. For each year, we compare the ex-post profit achieved when the electrolyzer does not participate in any ancillary service markets (denoted as ``none" in Figure \ref{fig:revenue_streams_all}), when it participates in either the mFRR or FCR capacity market (denoted as ``mFRR" or ``FCR"), and when it participates in both markets (denoted as ``both").

The ex-post profit increases when participating in any reserve market for all three years. When participating in both reserve markets, the electrolyzer increases its profits by 27\% in 2021 and 47\% in 2023. The most dramatic effect of reserve provision is observed in 2022 (middle plot). This is because 2022 was marked by high energy prices, and without the additional revenues from ancillary services, the electrolyzer would have incurred a financial loss while still being subject to the minimum hydrogen demand requirement. The increase in profit due to bidding in ancillary service markets can be attributed to different sources across the three years. In 2021, mFRR capacity prices were comparatively low, so the main increase in profit came from a combination of FCR capacity and reserve activation revenues. In 2022, FCR prices were extremely high, and the electrolyzer could not make a profit without bidding into the FCR capacity market. However, in 2023, FCR capacity prices declined, and very little additional revenue was generated from FCR service provision. Instead, a combination of mFRR capacity and activation payments provided an increase in profit.

The key takeaway is that while reserve provision increases the electrolyzer's profit in all three years, the additional revenue gained from the various markets differs significantly from year to year due to price volatility. We conclude that to reliably generate value through reserve provision, it is important for the electrolyzer to participate in multiple reserve markets.

\begin{figure}[!t]
    \centering
    \includegraphics[width=\linewidth]{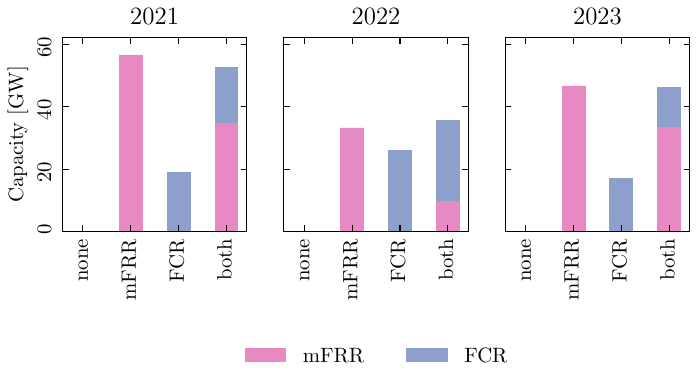}
    \caption{Total reserve capacity bids placed over the year, disregarding potential activations in the bidding stage, i.e., $\alpha = 0$.}
    \label{fig:capacities_three_years}
\end{figure}

\subsubsection{Reserve Capacity Bids Placed by the Electrolyzer}
\label{reserve-participation-section}

We present results on the total capacity bids (in GW) that the electrolyzer could place in the mFRR and FCR capacity markets over the years 2021-2023, while still disregarding potential activations at the bidding stage ($\alpha = 0$). These results are obtained directly from the optimization model \eqref{D2_obj}-\eqref{activation}, without requiring ex-post analysis.

Figure \ref{fig:capacities_three_years} illustrates the total capacity bids for the same four cases as in Figure \ref{fig:revenue_streams_all}. We observe that the total mFRR capacity bids remain significant even in 2021 and 2022, despite relatively low mFRR capacity prices during those years. Interestingly, the electrolyzer places more mFRR capacity bids in 2021 than in 2023, even though mFRR capacity prices were substantially higher in 2023. 
As discussed in Section \ref{analytical_results}, the electrolyzer may adjust its schedule to operate below its optimal setpoint in the day-ahead market to participate in the mFRR up reserve market, albeit incurring an opportunity cost. Consequently, during hours with low day-ahead prices, when the electrolyzer runs at full capacity, it can allocate its entire capacity to the mFRR reserve, even at very low reserve prices. The high frequency of low day-ahead price hours in 2021, as seen in Figure \ref{fig:prices}, and the asymmetric and hourly bid structure of the mFRR reserve therefore led to large capacities reserved for mFRR.

Our analytical results in Section \ref{analytical_results} also explain why the electrolyzer reserves a similar total quantity for the FCR market in 2021 and 2023, despite substantially higher reserve prices in 2021. This is because the higher day-ahead prices in 2021 increased the cost of participating in the reserve. Due to the four-hour block structure, even a single hour with sufficiently high day-ahead prices (e.g., 07:00-08:00 AM) could make FCR participation prohibitively expensive.

The key takeaway is that the current structure of mFRR and FCR capacity bids incentivizes an  electrolyzer exposed to day-ahead prices to participate in the mFRR market more frequently than in the FCR market, even when mFRR prices are lower, as observed in 2021 and 2022. 
\begin{newtext}
    This suggests that a relaxed FCR structure could increase the amount of electrolyzer capacity placed.
\end{newtext}
We explore this further in the following.

\begin{figure}[!t]
    \centering
    \includegraphics[width = 0.5\linewidth]{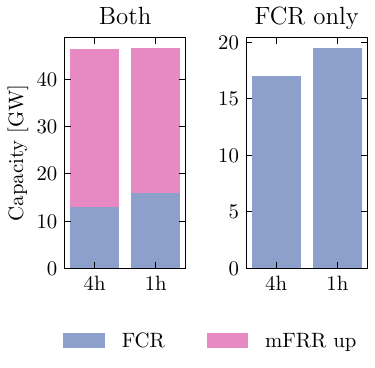}
    \caption{The total reserve capacity bids that the electrolyzer could place in 2023 under the current 4-hour FCR bid structure in DK1, and how these bids would change if the FCR bids were adjusted to an hourly structure, while still assuming $\alpha = 0$.}
    \label{fig:fcr_4h_1h}
\end{figure}

Figure \ref{fig:fcr_4h_1h} shows the total reserve capacity bids of the electrolyzer in 2023 for both the mFRR and FCR markets (left plot) and when bidding exclusively in the FCR market (right plot). In both cases, we compare two market structures for FCR bids: (\textit{i}) the current system with 4-hour blocks and (\textit{ii}) a hypothetical scenario allowing hourly FCR bids, referred to as the \textit{relaxed} FCR structure, as examined in Section \ref{analytical_results}. The results show that under this relaxed structure, the electrolyzer increases its FCR bidding by 23\% when participating in both markets and by 14\% when bidding only in the FCR market.

One may hypothesize that the electrolyzer would bid in the mFRR market less frequently if it considered potential activation at the bidding stage. This aspect will be analyzed further in Section \ref{activation-aware}.


\subsubsection{Impact of Activation on Hydrogen Demand}
\label{reserve-activation-section}

\begin{figure}[t]
    \centering
    \includegraphics[width=\linewidth]{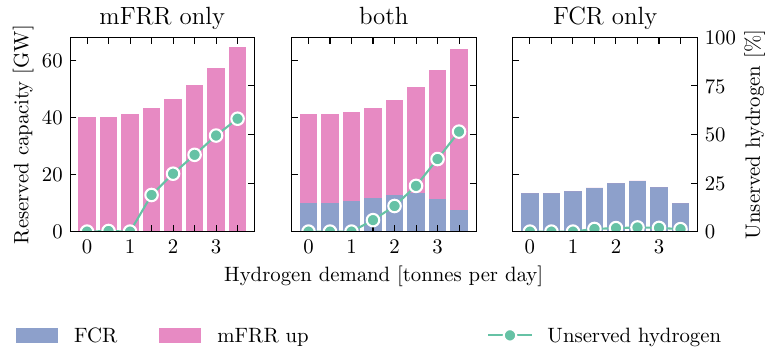}
    \caption{Total reserve capacity bids placed in 2023 under different minimum daily hydrogen demand levels, considering participation in mFRR only, FCR only, and both.}
    \label{fig:demand_sensitivity_alpha0}
\end{figure}

The scheduling results presented so far assume a daily minimum hydrogen demand of 2000 kg, which corresponds to slightly less than half of the electrolyzer’s daily production capacity. To better understand the impact of demand variations, we conduct a sensitivity analysis by adjusting the daily hydrogen demand in 500 kg increments, ranging from no minimum demand to 3500 kg per day. Additionally, we examine how reserve activation affects the ability to meet hydrogen demand. Due to data availability constraints, as explained in Section \ref{sec:histrotical_activation}, this analysis is limited to the year 2023.  The total capacity reserved in different markets under varying minimum daily hydrogen demand, assuming $\alpha = 0$, is presented in Figure \ref{fig:demand_sensitivity_alpha0} on the left axis. The green dots, corresponding to the right axis, shows the percentage of unserved minimum hydrogen demand.

Figure \ref{fig:demand_sensitivity_alpha0} shows that reducing the daily demand from the baseline of 2000 kg results in lower reserve capacity allocations for both products considered. At first glance, this may seem counterintuitive, as a lower hydrogen demand could imply greater operational flexibility for the electrolyzer. However, maintaining a minimum daily hydrogen output requires the electrolyzer to operate for a certain number of hours, regardless of day-ahead electricity prices. The scheduling model \eqref{D2_obj}-\eqref{activation} ensures that the expected hydrogen output meets or exceeds the minimum daily demand. Because both the mFRR up and the FCR reserves require the electrolyzer to be in operation to offer capacity, a lower daily demand reduces the number of hours with no opportunity cost associated with reserve participation when day-ahead prices are otherwise unfavorable for hydrogen production.  

Furthermore, Figure \ref{fig:demand_sensitivity_alpha0} shows that increasing the minimum daily demand above 2000 kg leads to a continued rise in mFRR capacity bids, whereas FCR participation declines. Once demand exceeds approximately half of the daily production capacity, the amount of reserves placed in FCR drops significantly. This highlights the inherent limitations of FCR in combination with hydrogen production due to its four-hour delivery blocks and symmetric bidding requirements. When the electrolyzer operates at (near) full capacity for extended periods, its ability to provide FCR capacity is constrained.  

These findings suggest that while stricter hydrogen demand requirements typically reduce operational flexibility, participation in hourly and asymmetrical reserve products like mFRR may increase under higher minimum demand scenarios.  

An increase in the minimum daily hydrogen demand also results in a higher percentage of unmet hydrogen demand. As illustrated in Figure \ref{fig:demand_sensitivity_alpha0}, the fraction of unserved hydrogen remains low ($<4$\%) when participating exclusively in the FCR market (right plot). This aligns with expectations based on our analysis in Section \ref{sec:histrotical_activation}. However, when participating in the mFRR market, unserved hydrogen demand becomes significant for minimum demand levels exceeding 1000 kg, ranging from approximately 10\% to 60\%.  

\begin{figure}[t]
    \centering
    \includegraphics[width=\linewidth]{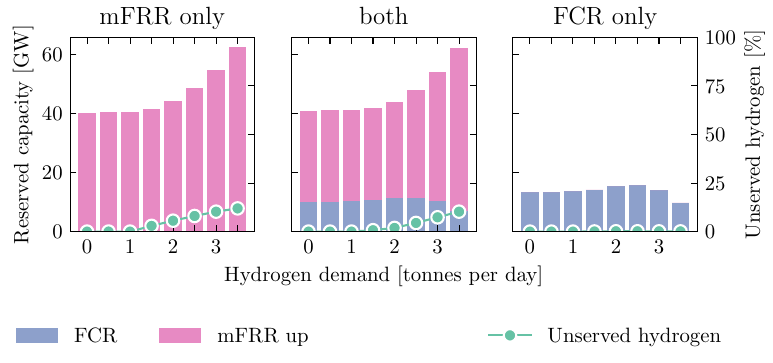}
    \caption{Reserves placed and unserved demand with weekly structure. The minimum demand levels are proportional to the daily demand, i.e., times seven across the week.}
    \label{fig:demand_relaxed}
\end{figure}

The extent to which reserve activation leads to unmet hydrogen demand will in reality depend on both the physical characteristics of the hydrogen supply infrastructure, such as the availability of a hydrogen pipeline, and contractual requirements, such as the need for continuous supply. These factors influence the definition of a \textit{demand period}. Figure \ref{fig:demand_relaxed} illustrates the capacity allocated to reserves and the resulting percentage of unserved hydrogen demand when considering a weekly demand period instead of a daily one.  Adopting a weekly minimum demand structure significantly reduces unmet hydrogen demand compared to a daily delivery period, with the shortfall decreasing to just 10\% at most. This result highlights the benefits of increased temporal flexibility on the hydrogen demand side for enabling reserve provision.

\subsection{Activation-Aware Reserve Participation ($\alpha > 0$)}
\label{activation-aware}

\begin{figure}[t]
    \centering
    \includegraphics[width=\linewidth]{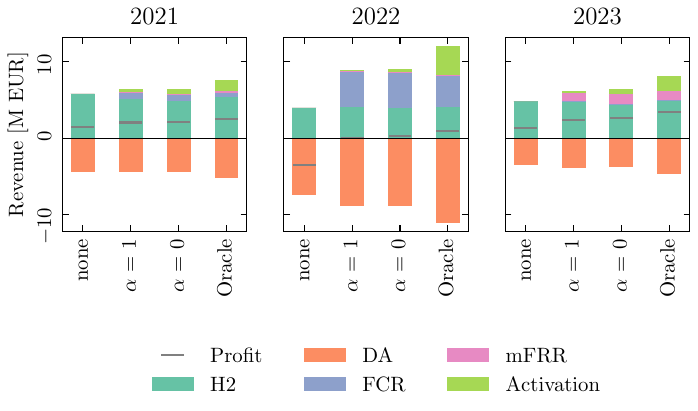}
    \caption{ \begin{newtext} Revenue streams and ex-post profits from participation in both reserve markets. For comparison, ``none" assumes no reserve market participation. In the cases of ``$\alpha = 1$" and ``$\alpha = 0$", the electrolyzer bids into both reserve markets, assuming full or no mFRR activation, respectively. The ``Oracle" case assumes perfect foresight of mFRR and FCR activations, as well as imbalance prices, and provides an upper bound on profit while satisfying the minimum hydrogen demand.
    \end{newtext}}
    \label{fig:profits_aware}
\end{figure}

\begin{figure}[t]
    \centering
    \includegraphics[width=\linewidth]{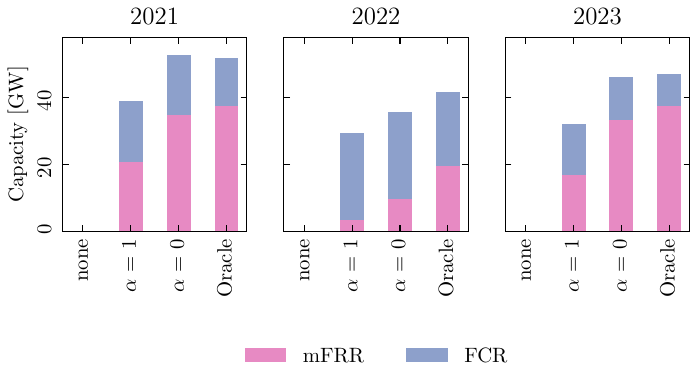}
    \caption{
    Total capacity bids placed when participating in both reserve markets.
    \begin{newtext}
         For comparison, ``none" assumes no reserve market participation. In the cases of ``$\alpha = 1$" and ``$\alpha = 0$", the electrolyzer bids into both reserve markets, assuming full or no mFRR activation, respectively. The ``Oracle" case assumes perfect foresight of mFRR and FCR activations, as well as imbalance prices, and provides an upper bound on profit while satisfying the minimum hydrogen demand.
    \end{newtext} }
    \label{fig:capacities_aware}
\end{figure}

As discussed in Section \ref{reserve-activation-section}, bidding capacity to the mFRR market without considering its potential activation can result in significant shares of unmet hydrogen demand (greater than 30\%) over a year, assuming a base demand of 2000 kg per day. To address this, we compare the results of the ``activation-unaware" model (i.e., $\alpha = 0$) to the ``activation-conservative" approach (i.e., $\alpha = 1$), as well as to the ``oracle" case, where the electrolyzer has perfect knowledge of future mFRR activation. The results are presented for all three years (2021-2023) under the minimum hydrogen demand of 2000 kg per day. Following the conclusion on the importance of participating in multiple markets, we focus on results for participation in both the FCR and mFRR markets. The revenue streams and ex-post profits are presented in Figure \ref{fig:profits_aware}, while the total reserve capacity bids placed are shown in Figure \ref{fig:capacities_aware}.

First, we examine the case of $\alpha = 1$, where capacity bids to the mFRR market are robust against activation. As shown in Figure \ref{fig:profits_aware}, participation in  reserve markets with activation-robust bidding still leads to increased profits compared to not participating in  reserve markets. This increase is primarily driven by FCR revenues in 2021 and 2022, while mFRR revenues contribute the most to the profit increase in 2023. Additionally, the profit increase under activation-robust bidding is comparable to the increase observed when activation is disregarded (i.e., $\alpha = 0$). Figure \ref{fig:capacities_aware} illustrates that, although activation-robust bidding reduces the total amount of capacity bid placed in the mFRR market, a significant mFRR capacity bid is still placed in 2021 and 2023. These years experienced generally lower day-ahead prices, meaning the electrolyzer often schedules consumption above the minimum daily hydrogen demand. As a result, the electrolyzer can still bid its excess capacity to the mFRR market while remaining robust against potential activation impacts on minimum daily hydrogen demand. Interestingly, 2022 stands out as the year with both the lowest mFRR revenue (under all $\alpha$ assumptions) and the largest profit difference between $\alpha = 1$ and $\alpha = 0$. The combination of low mFRR prices and high day-ahead prices in 2022 led to minimal mFRR capacity bids in the activation-robust case. However, the activation of mFRR services generated substantial revenues due to high balancing prices, leading to a larger profit difference between the $\alpha = 1$ and $\alpha = 0$ cases.

Next, we analyze the ``oracle" case, where the electrolyzer has perfect information about future mFRR activations and balancing prices at the time of bidding. Thus, the electrolyzer can place bids that maximize profits without violating the minimum daily demand constraint. Note that this constraint is still strictly enforced, with no flexibility for missed demand through penalties or other means. The ``oracle" leads to a significant increase in profit throughout all years compared to the other cases, as seen in Figure \ref{fig:profits_aware}. This increase in profit results from capacity placed in the mFRR market during hours of high balancing prices, which generates large revenues upon activation. 

In conclusion, the electrolyzer can increase its profits through ancillary service provision while maintaining robustness against the effects of reserve activation on its hydrogen demand obligations.


\section{Conclusion} \label{conclusion}

This work illustrates how electrolyzers can increase their operational profits through secondary revenue streams from ancillary service provision. We specifically investigate two different reserve products, the fast and low-energy FCR reserve, and the slower, energy-intensive mFRR reserve, on a case study formed by historical data from 2021-2023 in the western Danish bidding zone DK1.

First, we show analytically that the electrolyzer's opportunity cost of offering a reserve, i.e., its reserve bid price, is a function of the day-ahead price under certain assumptions. Further, we show the effect of the reserve structure on the electrolyzer's bid price, which explains why the mFRR is more suited for electrolyzer participation than the FCR. Second, we show how participation in reserves can significantly improve the operational profits of the electrolyzer, even when the electrolyzer participates in an activation-robust manner. Finally, we illustrate the impact of reserve activation on meeting a minimum hydrogen demand, showing the benefit of temporal flexibility in the demand structure when the electrolyzer participates in reserve markets.  


We make two recommendations to relevant stakeholders. The TSOs should consider the reserve bid structures if they expect large shares of reserves to be provided by electrolyzers, as the reserve bid structure might have a significant impact on how much capacity electrolyzers will offer to the market and at what price. Electrolyzer owners should not necessarily disregard energy-intensive reserves such as mFRR as it is possible to make a profit while ensuring that the hydrogen demand is met upon activation, and consider the benefits of relaxed hydrogen demand structures when negotiating off-take contracts.





Future work should explore the assumption of perfect price foresight and assess the extent to which the additional value of reserve provision might be influenced by uncertainty in prices. This is particularly important, as we have shown that the reserve bid price of the electrolyzer depends on the (forecasted) day-ahead price. Additionally, the long-term impacts of reserve participation on electrolyzer performance and maintenance costs due to, e.g., degradation of the electrolyzer cells and potential reductions in lifetime, should be evaluated. However, obtaining accurate data to model these aspects has proven challenging for the authors. Finally, an analysis of price cannibalization as more technologies enter the reserve markets should be conducted.

\appendix

\section{FCR Activation}
\label{app:FCR}
The \textit{real} FCR activation signal at second $\tau$, denoted $\alpha_{\tau}^{\rm{real, F}}$,  is a function of the frequency measurement $f_{\tau}$ and is calculated as follows:
\begin{align}
     \alpha^{\rm{real, F}}_{\tau} = 
\begin{cases}
    1, & \text{if } \ f_{\tau} \geq 50.1  \\
    - 1, & \text{if } \ f_{\tau} \leq 49.9  \\
    \frac{f_{\tau} - 50.0} {0.1},              & \text{otherwise.}
\end{cases}
\end{align}

Therefore, the  change in electrolyzer's power consumption, $\Delta P^{\rm{F}}_t$, and the corresponding change in energy consumption, $\Delta E^{\rm{F}}_t$, at hour $t$ are computed as follows:  
\begin{subequations} 
\begin{align}
    &\Delta P^{\rm{F}}_{t} = \frac{r^{\rm{F*}}_{t}}{3600} \sum_{\tau  \in  T_{t}}\alpha^{\rm{real, F}}_{\tau}, \ \ \ \forall  t \in \mathcal{T}\\
    &\Delta E^{\rm{F}}_{t} = \Delta P^{\rm{F}}_{t} \Delta_{\rm{T}}, \ \ \ \forall  t \in \mathcal{T},
\end{align}
\end{subequations} 
where $\tau \in T_t$ denotes the set of 3600 seconds within hour $t$, with a time step of $\Delta_{\rm{T}} = 1$. The FCR activation during hour $t$ adjusts the electrolyzer's power consumption as a fraction of its optimal FCR quantity bid, $r^{\rm{F*}}_t$, determined from \eqref{D2_obj}-\eqref{activation}, which has already been placed in the FCR market. When the grid frequency exceeds 50 Hz ($\alpha_{\tau}^{\rm{real, F}} > 0$), the electrolyzer is automatically activated to consume additional electricity, resulting in increased hydrogen production. Conversely, when the frequency drops below 50 Hz ($\alpha_{\tau}^{\rm{real, F}} < 0$), the electrolyzer's power consumption decreases, leading to reduced hydrogen production.

\section{mFRR Activation}
\label{app:mFRR}
The Danish transmission system operator, Energinet, is currently procuring only upward mFRR reserves, i.e., reserves for addressing a supply deficit in the system. Let $I_t$ denote the system imbalance at hour $t$, where $I_t < 0$ indicates the system's need for up-regulation. The \textit{realized} upward activation signal, $\alpha^{\rm{real, m \uparrow}}$, at hour $t$ is then computed as:

\begin{align}    
     \alpha^{\rm{real, m \uparrow}}_{t} = 
\begin{cases}
    1, & \text{if } \ \Bar{\eta} \lambda^{\rm{h}} \leq \lambda^{\rm{B}}_{t} \text{ and } I_t < 0  \\
    0,              & \text{otherwise,}
\end{cases}
\end{align}
where $\alpha^{\rm{real, m \uparrow}}_{t} > 0$ at hour $t$ refers to upward activation, implying a decrease in the power consumption of the electrolyzer. Recall that for FCR activation, it was the opposite, meaning $\alpha_{\tau}^{\rm{real, F}} > 0$ at second $\tau$ implied an increase in power consumption.

The change in the electrolyzer's power consumption as a consequence of upward mFRR activation, $\Delta P^{\rm{m \uparrow}}$, and the subsequent change in energy consumption, $\Delta E^{\rm{m \uparrow}}$, are then computed as:
\begin{subequations} 
\begin{align}
    &\Delta P^{\rm{m \uparrow }}_{t} = -r^{\rm{m \uparrow*}}_{t} \alpha^{\rm{real, m \uparrow}}_{t},  \ \ \ \forall  t \in \mathcal{T} \label{cd}\\
    &\Delta E^{\rm{m \uparrow }}_{t} = \Delta P^{\rm{m \uparrow }}_{t} \Delta_{\rm{T}},  \ \ \ \forall  t \in \mathcal{T},
\end{align}
\end{subequations} 
where $r^{\rm{m \uparrow*}}_{t}$ represents the optimal upward mFRR bid placed at hour $t$, determined from \eqref{D2_obj}-\eqref{activation}. Note that, with the minus sign included in \eqref{cd}, an upward mFRR activation ($\alpha^{\rm{real, m \uparrow}}_{t} > 0$) results in reduced power consumption ($\Delta P^{\rm{m \uparrow }}_{t} <0$) and, consequently, decreased hydrogen production ($\Delta E^{\rm{m \uparrow }}_{t} <0$).

\section{Unmet Hydrogen Calculation}
\label{app:unmet_h2}
If the electrolyzer reduces its power consumption in response to an activation event 
($\Delta p_{t}^{(.)} < 0$), it will also reduce its hydrogen production ($\Delta h_{t}^{(.)} < 0$), and vice versa.  The final term in \eqref{eq:profit_ex_post} represents a penalty, $h^{\rm{unmet}}$, which arises if the electrolyzer fails to meet the minimum hydrogen demand, \underline{H}, as specified in the bilateral purchase agreement with off-takers. Note that the proposed optimization model \eqref{D2_obj}-\eqref{activation} ensures that this minimum hydrogen requirement is fully met. Therefore, any unmet hydrogen obligation resulting from reserve activation is penalized.  The unmet hydrogen obligation is computed as follows:
\begin{align}
    h^{\rm{unmet}} = {\rm max}\Bigg \{  \underline{H} - \sum_{t \in \mathcal{T}}{\rm H}\Big(p^{\rm tot}_t  + \Delta P^{\rm{F}}_{t} + \Delta P^{\rm{m \uparrow }}_{t}\Big), 0 \Bigg \},
\end{align}
where any positive value of $h^{\rm{unmet}}$ indicates a failure to meet the minimum hydrogen demand requirement $\underline{H}$. 

{\footnotesize
\nomenclature[V]{\( \)}{{$^*$most variables are given per hour $t$  }
\nomunit{}}

\nomenclature[V]{\(h^{\rm act \uparrow}_t\)}{{Hydrogen production upon up-reserve activation}
\nomunit{[kg]}}

\nomenclature[V]{\(h^{\rm act \downarrow}_t\)}{{Hydrogen production upon down-reserve activation }
\nomunit{[kg]}}

\nomenclature[V]{\(h^{\rm d}_{t,d}\)}{{Hydrogen dispensed to tube-trailer $d$}
\nomunit{[kg]}}

\nomenclature[V]{\(h^{\rm e}_t\)}{{Expected hydrogen production (pre-activation)}
\nomunit{[kg]}}

\nomenclature[V]{\(p^{\rm act \uparrow}_{t,s}\)}{{Electrolyzer consumption upon up activation at segment $s$}
\nomunit{[MW]}}

\nomenclature[V]{\(p^{\rm act \downarrow}_{t,s}\)}{{Electrolyzer consumption upon down activation at segment $s$}
\nomunit{[MW]}}

\nomenclature[V]{\(p^{\rm DA}_t\)}{{Electrolyzer day-ahead power bid }
\nomunit{[MW]}}

\nomenclature[V]{\(p^{\rm e}_{t,s}\)}{{Electrolyzer consumption at segment $s$}
\nomunit{[MW]}}

\nomenclature[V]{\(r^{\rm F}_i\)}{{FCR bid of block $i$}
\nomunit{[MW]}}

\nomenclature[V]{\(r^{\rm m \uparrow}_t\)}{{mFRR up bid }
\nomunit{[MW]}}

\nomenclature[V]{\(r^{\rm m \downarrow}_t\)}{{mFRR down bid }
\nomunit{[MW]}}

\nomenclature[V]{\(p^{\rm sb}_t\)}{{Electrolyzer stand-by consumption}
\nomunit{[MW]}}

\nomenclature[V]{\(p^{\rm tot}_t\)}{{Electrolyzer total consumption}
\nomunit{[MW]}}

\nomenclature[V]{\(s^{\rm d}_{t,d}\)}{{Storage state at dispenser $d$}
\nomunit{[MW]}}

\nomenclature[V]{\(z^{\rm act \uparrow}_{t,s}\)}{{Active segment $s$ indicator upon up-activation}
\nomunit{\{0,1\}}}

\nomenclature[V]{\(z^{\rm act \downarrow}_{t,s}\)}{{Active segment $s$ indicator upon down-activation}
\nomunit{\{0,1\}}}

\nomenclature[V]{\(z^{\rm e}_{t,s}\)}{{Active segment $s$ indicator}
\nomunit{\{0,1\}}}

\nomenclature[V]{\(z^{\rm F}_i\)}{{FCR bid indicator in block $i$}
\nomunit{\{0,1\}}}

\nomenclature[V]{\(z^{\rm m \uparrow}_t\)}{{mFRR-up bid indicator}
\nomunit{\{0,1\}}}

\nomenclature[V]{\(z^{\rm m \downarrow}_t\)}{{mFRR-down bid indicator}
\nomunit{\{0,1\}}}

\nomenclature[V]{\(z^{\rm off}_t\)}{{Electrolyzer off-state indicator}
\nomunit{\{0,1\}}}

\nomenclature[V]{\(z^{\rm on}_t\)}{{Electrolyzer on-state indicator}
\nomunit{\{0,1\}}}

\nomenclature[V]{\(z^{\rm sb}_t\)}{{Electrolyzer standby-state indicator}
\nomunit{\{0,1\}}}

\nomenclature[P]{\(\lambda^{\rm h}\)}{{Contracted hydrogen price}
\nomunit{[€/kg]}}

\nomenclature[P]{\(\lambda^{\rm m\uparrow}_t\)}{{mFRR-up price in hour $t$}
\nomunit{[€/MW/h]}}

\nomenclature[P]{\(\lambda^{\rm m\downarrow}_t\)}{{mFRR-down price in hour $t$}
\nomunit{[€/MW/h]}}

\nomenclature[P]{\(\lambda^{\rm DA}_{t}\)}{{Day-ahead price in hour $t$}
\nomunit{[€/MWh]}}

\nomenclature[P]{\(\lambda^{\rm F}_i\)}{{FCR price in block $i$}
\nomunit{[€/MW/h]}}

\nomenclature[P]{\(\Delta T^{\rm F}\)}{{FCR bid-block duration}
\nomunit{[h]}}

\nomenclature[P]{\(A_s\)}{{Slope of segment $s$}
\nomunit{[kg/MW]}}

\nomenclature[P]{\(B_s\)}{{Intercept of segment $s$}
\nomunit{[kg]}}

\nomenclature[P]{\(P^{\rm sb}\)}{{Electrolyzer standby consumption}
\nomunit{[MW]}}

\nomenclature[P]{\(C^{e}\)}{{Electrolyzer capacity}
\nomunit{[MW]}}

\nomenclature[P]{\(\underline{E}\)}{{Electrolyzer minimum consumption}
\nomunit{[MW]}}

\nomenclature[P]{\(\underline{P}^{\rm F}\)}{{FCR minimum  bid size}
\nomunit{[MW]}}

\nomenclature[P]{\(\underline{P}^{\rm m}\)}{{mFRR minimum  bid size}
\nomunit{[MW]}}

\nomenclature[P]{\(\overline{P}^{\rm F}\)}{{FCR maximum  bid size}
\nomunit{[MW]}}

\nomenclature[P]{\(\overline{P}^{\rm m}\)}{{mFRR maximum  bid size}
\nomunit{[MW]}}

\nomenclature[P]{\(\overline{P}_s\)}{{Maximum power at segment $s$}
\nomunit{[MW]}}

\nomenclature[P]{\(\underline{P}_s\)}{{Minimum power at segment $s$}
\nomunit{[MW]}}

\nomenclature[P]{\(\underline{H}\)}{{Minimum daily hydrogen demand}
\nomunit{[kg]}}

\nomenclature[P]{\(C^{\rm d}\)}{{Dispenser capacity}
\nomunit{[kg/h]}}

\nomenclature[P]{\(S^{\rm d}\)}{{Storage (tube trailer) capacity}
\nomunit{[kg]}}

\nomenclature[P]{\(\zeta_{t,d}
\)}{{Tube trailer availability at dispenser $d$ in hour $t$}
\nomunit{\{0,1\}}}

\nomenclature[P]{\(\alpha^{\rm{act \uparrow}}_{t}
\)}{{Assumed upward activation in hour $t$}
\nomunit{[\%]}}

\nomenclature[P]{\(\alpha^{\rm{act \downarrow}}_{t}
\)}{{Assumed downward activation in hour $t$}
\nomunit{[\%]}}


\nomenclature[P]{\(\tilde{p}^{\rm DA*}_{t}\)}{{Optimal day-ahead schedule without reserve participation}
\nomunit{[MW]}}

\nomenclature[F]{\({\rm p}^{\rm \downarrow}(\cdot) \)}
{{New setpoint as a function of down reserve capacity}
\nomunit{[MW]}}

\nomenclature[F]{\({\rm p}^{\rm \uparrow}(\cdot)\)}
{{New setpoint as a function of up reserve capacity}
\nomunit{[MW]}}

\nomenclature[F]{\({\rm H}(\cdot)
\)}{{Empirical hydrogen production as a function of a setpoint}
\nomunit{[kg/h]}}

\nomenclature[F]{\({\rm h}(\cdot)
\)}{{Linearized hydrogen production as a function of a setpoint}
\nomunit{[kg/h]}}

\nomenclature[F]{\(\tilde{\Gamma}(\cdot)\)}{{ Profit as a function of a setpoint (without reserve revenues) }
\nomunit{[EUR]}}

\nomenclature[F]{\({\rm C}^{\downarrow}(\cdot)\)}{{ Opportunity cost as a function of a down reserve capacity}
\nomunit{[EUR]}}

\nomenclature[F]{\({\rm C}^{\uparrow}(\cdot)\)}{{ Opportunity cost as a function of a up reserve capacity}
\nomunit{[EUR]}}

\nomenclature[F]{\(\pi^{\downarrow}(\cdot)
\)}{{ Bid price as a function of a down reserve capacity }
\nomunit{[EUR/MW]}}

\nomenclature[F]{\(\pi^{\uparrow}(\cdot)
\)}{{ Bid price as a function of a up reserve capacity }
\nomunit{[EUR/MW]}}

\nomenclature[F]{\(\pi^{\downarrow\uparrow}(\cdot)
\)}{{ Bid price as a function of a symmetric reserve capacity }
\nomunit{[EUR/MW]}}

\nomenclature[F]{\(\pi^{\rm \downarrow\uparrow, 4h}(\cdot)
\)}{{ Bid price as a function of a 4-hour, symmetric reserve }
\nomunit{[EUR/MW]}}

\nomenclature[S]{\(\mathcal{I} \)}{{Set of FCR reserve blocks }
\nomunit{}}

\nomenclature[S]{\(\mathcal{S} \)}{{Set of segments on linearized hydrogen production curve }
\nomunit{}}

\nomenclature[S]{\(\mathcal{T}_i \)}{{Set of hours in reserve block $i$ }
\nomunit{}}

\nomenclature[S]{\(\mathcal{T} \)}{{Set of hours in a day }
\nomunit{}}

\nomenclature[S]{\(\mathcal{D} \)}{{Set of hydrogen dispensers }
\nomunit{}}

\nomenclature[S]{
\begin{newtext}
\(\mathcal{N}_t \)
\end{newtext}
}{{
\begin{newtext}
Set of $N^{\rm down}$ hours from hour $t$ (electrolyzer down-time)
\end{newtext}
}
\nomunit{}}

\nomenclature[V]{\(Z\)}{{ Objective value}
\nomunit{[EUR]}}

\nomenclature[P]{\(\Gamma\)}{{ Ex-post profit upon historical activation}
\nomunit{[EUR]}}

\nomenclature[P]{\(\lambda^{\rm{B}}_{t}
\)}{{Balancing price in hour $t$}
\nomunit{[EUR/MWh]}}

\nomenclature[P]{\(\overline\eta
\)}{{Peak electrolyzer efficiency}
\nomunit{[kg/h/MW]}}

\nomenclature[P]{\(\Delta E^{\rm F}_t
\)}{{Change in energy consumption due to FCR activation in hour $t$}
\nomunit{[MWh]}}

\nomenclature[P]{\(\Delta E^{\rm m\uparrow}_t
\)}{{Change in energy consumption due to mFRR up activation in hour $t$}
\nomunit{[MWh]}}

\nomenclature[P]{\(\Delta P^{\rm F}_t
\)}{{Change in power setpoint due to FCR activation in hour $t$}
\nomunit{[MW]}}

\nomenclature[P]{\(\Delta P^{\rm m\uparrow}_t
\)}{{Change in power setpoint due to mFRR up activation in hour $t$}
\nomunit{[MW]}}

\nomenclature[P]{\(\Delta h_t\)}{{Change in hydrogen production due to reserve activation in hour $t$}
\nomunit{[kg]}}

\printnomenclature

}
\vspace{2mm}
\section*{Acknowledgement }
We thank the Danish Energy Technology Development and Demonstration Programme for supporting this research through HOMEY project (Grant number: 64021-7010). 
\vspace{2mm}



\Urlmuskip=0mu plus 1mu\relax

\bibliography{references.bib}
\bibliographystyle{elsarticle-harv} 




\end{document}